\documentclass[11pt]{article}
\usepackage[dvips]{graphicx}
\usepackage{textcomp}
\usepackage{amsbsy}
\usepackage{latexsym}
\usepackage[mathscr]{eucal}
\usepackage{amsfonts,amsmath,amsthm}
\usepackage{amssymb}
\usepackage[usenames]{xcolor}

\usepackage{fullpage}

\begin{document}
\bibliographystyle{plain}
\title
{Lipschitz widths}
\author{ 
Guergana Petrova,  Przemys{\l}aw Wojtaszczyk
\thanks{%
   G.P.  was supported by the  NSF Grant  DMS 2134077, Tripods Grant CCF-1934904, and ONR Contract N00014-20-1-278. P. W. was supported by National Science Centre, 
   Polish grant UMO-2016/21/B/ST1/00241.  }}
\hbadness=10000
\vbadness=10000
\newtheorem{lemma}{Lemma}[section]
\newtheorem{prop}[lemma]{Proposition}
\newtheorem{cor}[lemma]{Corollary}
\newtheorem{theorem}[lemma]{Theorem}
\newtheorem{remark}[lemma]{Remark}
\newtheorem{example}{Example}[section]
\newtheorem{definition}[lemma]{Definition}
\newtheorem{proper}[lemma]{Properties}
\newtheorem{assumption}[lemma]{Assumption}
%
\newenvironment{disarray}{\everymath{\displaystyle\everymath{}}\array}{\endarray}

\def\RR{\rm \hbox{I\kern-.2em\hbox{R}}}
\def\NN{\rm \hbox{I\kern-.2em\hbox{N}}}
\def\ZZ{\rm {{\rm Z}\kern-.28em{\rm Z}}}
\def\CC{\rm \hbox{C\kern -.5em {\raise .32ex \hbox{$\scriptscriptstyle
|$}}\kern
-.22em{\raise .6ex \hbox{$\scriptscriptstyle |$}}\kern .4em}}
\def\vp{\varphi}
\def\<{\langle}
\def\>{\rangle}
\def\t{\tilde}
\def\i{\infty}
\def\e{\varepsilon}
\def\sm{\setminus}
\def\nl{\newline}
\def\o{\overline}
\def\wt{\widetilde}
\def\wh{\widehat}
\def\cT{{\cal T}}
\def\cA{{\cal A}}
\def\cI{{\cal I}}
\def\cV{{\cal V}}
\def\cB{{\cal B}}
\def\cF{{\cal F}}
\def\cY{{\cal Y}}

\def\cD{{\cal D}}
\def\cP{{\cal P}}
\def\cJ{{\cal J}}
\def\cM{{\cal M}}
\def\cO{{\cal O}}
\def\Chi{\raise .3ex
\hbox{\large $\chi$}} \def\vp{\varphi}
\def\lsima{\hbox{\kern -.6em\raisebox{-1ex}{$~\stackrel{\textstyle<}{\sim}~$}}\kern -.4em}
\def\lsim{\hbox{\kern -.2em\raisebox{-1ex}{$~\stackrel{\textstyle<}{\sim}~$}}\kern -.2em}
\def\[{\Bigl [}
\def\]{\Bigr ]}
\def\({\Bigl (}
\def\){\Bigr )}
\def\[{\Bigl [}
\def\]{\Bigr ]}
\def\({\Bigl (}
\def\){\Bigr )}
\def\L{\pounds}
\def\pr{{\rm Prob}}
\newcommand{\cs}[1]{{\color{magenta}{#1}}}
\def\ds{\displaystyle}
\def\ev#1{\vec{#1}}     
\newcommand{\lt}{\ell^{2}(\nabla)}
\def\Supp#1{{\rm supp\,}{#1}}
\def\R{\mathbb{R}}
\def\E{\mathbb{E}}
\def\nl{\newline}
\def\T{{\relax\ifmmode I\!\!\hspace{-1pt}T\else$I\!\!\hspace{-1pt}T$\fi}}
\def\N{\mathbb{N}}
\def\Z{\mathbb{Z}}
\def\N{\mathbb{N}}
\def\Zd{\Z^d}
\def\Q{\mathbb{Q}}
\def\C{\mathbb{C}}
\def\Rd{\R^d}
\def\gsim{\mathrel{\raisebox{-4pt}{$\stackrel{\textstyle>}{\sim}$}}}
\def\sime{\raisebox{0ex}{$~\stackrel{\textstyle\sim}{=}~$}}
\def\lsim{\raisebox{-1ex}{$~\stackrel{\textstyle<}{\sim}~$}}
\def\div{\mbox{ div }}
\def\M{M}  \def\NN{N}                  
\def\Le{{\ell^1}}            
\def\Lz{{\ell^2}}
\def\Let{{\tilde\ell^1}}     
\def\Lzt{{\tilde\ell^2}}
\def\Ltw{\ell^\tau^w(\nabla)}
\def\t#1{\tilde{#1}}
\def\la{\lambda}
\def\La{\Lambda}
\def\ga{\gamma}
\def\BV{{\rm BV}}
\def\Ga{\eta}
\def\al{\alpha}
\def\cZ{{\cal Z}}
\def\cA{{\cal A}}
\def\cU{{\cal U}}
\def\argmin{\mathop{\rm argmin}}
\def\argmax{\mathop{\rm argmax}}
\def\prob{\mathop{\rm prob}}

\def\cO{{\cal O}}
\def\cA{{\cal A}}
\def\cC{{\cal C}}
\def\cS{{\cal F}}
\def\bu{{\bf u}}
\def\bz{{\bf z}}
\def\bZ{{\bf Z}}
\def\bI{{\bf I}}
\def\cE{{\cal E}}
\def\cD{{\cal D}}
\def\cG{{\cal G}}
\def\cI{{\cal I}}
\def\cJ{{\cal J}}
\def\cM{{\cal M}}
\def\cN{{\cal N}}
\def\cT{{\cal T}}
\def\cU{{\cal U}}
\def\cV{{\cal V}}
\def\cW{{\cal W}}
\def\cL{{\cal L}}
\def\cB{{\cal B}}
\def\cG{{\cal G}}
\def\cK{{\cal K}}
\def\cX{{\cal X}}
\def\cS{{\cal S}}
\def\cP{{\cal P}}
\def\cQ{{\cal Q}}
\def\cR{{\cal R}}
\def\cU{{\cal U}}
\def\bL{{\bf L}}
\def\bl{{\bf l}}
\def\bK{{\bf K}}
\def\bC{{\bf C}}
\def\X{X\in\{L,R\}}
\def\ph{{\varphi}}
\def\D{{\Delta}}
\def\H{{\cal H}}
\def\bM{{\bf M}}
\def\bx{{\bf x}}
\def\bj{{\bf j}}
\def\bG{{\bf G}}
\def\bP{{\bf P}}
\def\bW{{\bf W}}
\def\bT{{\bf T}}
\def\bV{{\bf V}}
\def\bv{{\bf v}}
\def\bt{{\bf t}}
\def\bz{{\bf z}}
\def\bw{{\bf w}}
\def \span{{\rm span}}
\def \meas {{\rm meas}}
\def\rhom{{\rho^m}}
\def\diff{\hbox{\tiny $\Delta$}}
\def\EE{{\rm Exp}}
\def\lll{\langle}
\def\argmin{\mathop{\rm argmin}}
\def\codim{\mathop{\rm codim}}
\def\rank{\mathop{\rm rank}}

\def\argmax{\mathop{\rm argmax}}
\def\dJ{\nabla}
\newcommand{\ba}{{\bf a}}
\newcommand{\bb}{{\bf b}}
\newcommand{\bc}{{\bf c}}
\newcommand{\bd}{{\bf d}}
\newcommand{\bs}{{\bf s}}
\newcommand{\bff}{{\bf f}}
\newcommand{\bp}{{\bf p}}
\newcommand{\bg}{{\bf g}}
\newcommand{\by}{{\bf y}}
\newcommand{\br}{{\bf r}}
\newcommand{\be}{\begin{equation}}
\newcommand{\ee}{\end{equation}}
\newcommand{\bea}{$$ \begin{array}{lll}}
\newcommand{\eea}{\end{array} $$}
\def \Vol{\mathop{\rm  Vol}}
\def \mes{\mathop{\rm mes}}
\def \Prob{\mathop{\rm  Prob}}
\def \exp{\mathop{\rm    exp}}
\def \sign{\mathop{\rm   sign}}
\def \sp{\mathop{\rm   span}}
\def \rad{\mathop{\rm   rad}}
\def \vphi{{\varphi}}
\def \csp{\overline \mathop{\rm   span}}
%
%
\newcommand{\beqn}{\begin{equation}}
\newcommand{\eeqn}{\end{equation}}
\def\beginproof{\noindent{\bf Proof:}~ }
\def\endproof{\hfill\rule{1.5mm}{1.5mm}\\[2mm]}

\newenvironment{Proof}{\noindent{\bf Proof:}\quad}{\endproof}

\renewcommand{\theequation}{\thesection.\arabic{equation}}
\renewcommand{\thefigure}{\thesection.\arabic{figure}}

\makeatletter
\@addtoreset{equation}{section}
\makeatother

\newcommand\abs[1]{\left|#1\right|}
\newcommand\clos{\mathop{\rm clos}\nolimits}
\newcommand\trunc{\mathop{\rm trunc}\nolimits}
\renewcommand\d{d}
\newcommand\dd{d}
\newcommand\diag{\mathop{\rm diag}}
\newcommand\dist{\mathop{\rm dist}}
\newcommand\diam{\mathop{\rm diam}}
\newcommand\cond{\mathop{\rm cond}\nolimits}
\newcommand\eref[1]{{\rm (\ref{#1})}}
\newcommand{\iref}[1]{{\rm (\ref{#1})}}
\newcommand\Hnorm[1]{\norm{#1}_{H^s([0,1])}}
\def\int{\intop\limits}
\renewcommand\labelenumi{(\roman{enumi})}
\newcommand\lnorm[1]{\norm{#1}_{\ell^2(\Z)}}
\newcommand\Lnorm[1]{\norm{#1}_{L_2([0,1])}}
\newcommand\LR{{L_2(\R)}}
\newcommand\LRnorm[1]{\norm{#1}_\LR}
\newcommand\Matrix[2]{\hphantom{#1}_#2#1}
\newcommand\norm[1]{\left\|#1\right\|}
\newcommand\ogauss[1]{\left\lceil#1\right\rceil}
\newcommand{\QED}{\hfill
\raisebox{-2pt}{\rule{5.6pt}{8pt}\rule{4pt}{0pt}}%
  \smallskip\par}
\newcommand\Rscalar[1]{\scalar{#1}_\R}
\newcommand\scalar[1]{\left(#1\right)}
\newcommand\Scalar[1]{\scalar{#1}_{[0,1]}}
\newcommand\Span{\mathop{\rm span}}
\newcommand\supp{\mathop{\rm supp}}
\newcommand\ugauss[1]{\left\lfloor#1\right\rfloor}
\newcommand\with{\, : \,}
\newcommand\Null{{\bf 0}}
\newcommand\bA{{\bf A}}
\newcommand\bB{{\bf B}}
\newcommand\bR{{\bf R}}
\newcommand\bD{{\bf D}}
\newcommand\bE{{\bf E}}
\newcommand\bF{{\bf F}}
\newcommand\bH{{\bf H}}
\newcommand\bU{{\bf U}}
\newcommand \A {{\bb A}}
\newcommand\cH{{\cal H}}
\newcommand\sinc{{\rm sinc}}
\def\enorm#1{| \! | \! | #1 | \! | \! |}

\newcommand{\am}{a_{\min}}
\newcommand{\aM}{a_{\max}}

\newcommand{\dm}{\frac{d-1}{d}}

\let\bm\bf
\newcommand{\bbeta}{{\mbox{\boldmath$\beta$}}}
\newcommand{\bal}{{\mbox{\boldmath$\alpha$}}}
\newcommand{\bbi}{{\bm i}}

\def\nnew{\color{Red}}
\def\onew{\color{orange}}
\def\wnew{\color{magenta}}

\newcommand{\dI}{\Delta}
\newcommand\aconv{\mathop{\rm absconv}}

\maketitle
 %

\vskip 0.2in

\centerline{\it To Ron DeVore, with the utmost respect and admiration}

  \begin{abstract}   
  
  This paper  introduces a  measure,  called  Lipschitz widths, of the optimal performance possible of certain nonlinear methods of approximation. It discusses their relation to  entropy numbers and other well known widths such as the Kolmogorov and the stable manifold widths. It also shows that the Lipschitz widths provide a theoretical benchmark for the approximation quality achieved via deep neural networks.
  
  
        \end{abstract}

\section{Introduction}

Nonlinear methods of approximation provide reliable and efficient ways of investigating  the underlying phenomena in many application areas. Despite of their extensive usage however, there is still a lack of comprehensive  understanding of the  intrinsic limitations of these nonlinear methods, even on a  purely  theoretical level. Several mathematical concepts, called widths, have been  established to access numerous aspects of the quality of linear and nonlinear  approximations. 
As such, we mention the classical by now Kolmogorov, linear, manifold, Gelfand widths, which give a theoretical benchmark on what is the best possible performance of particular  methods of approximation, see
 \cite{DKLT}, where  a summary of different nonlinear widths and their relations to one another is discussed.  
  
  Recently,   Deep Neural Networks (DNN) have been used extensively as a method  of  choice  for variety of machine learning problems  and  as a computational  platform in many other areas. Despite of their empirical  successes, the explanation of the reasons behind their stellar performance is still in its infancy.  On mathematical level,  DNN   can be viewed as a method of nonlinear approximation of an  underlying function $f$, where the approximant $\Phi(y)\approx f$ is  a continuous function, generated by a DNN with parameters $y$. It can be shown that the mapping which to every choice of parameters $y$ of the DNN assigns  $\Phi(y)$ is in fact a Lipschitz mapping. Thus, DNN approximation is a particular case of a nonlinear approximation of a function $f$, or a compact class $\cK$, by  the images  of Lipschitz mappings. Then,  the question of DNN optimal performance  is intimately related  to the quantification  of the optimal performance of such nonlinear methods and to the introduction and study of corresponding ways to measure it. 
 A width, called {\it stable manifold}  width,  was  presented in \cite{CDPW}, with the sole purpose to determine   the optimal performance of such nonlinear methods in the context of numerical computation, where the stability plays an essential role.  In this paper,  we take a slightly different point of view and introduce  the concept of {\it Lipschitz} widths, where we are not so concerned about  the numerical stability of the method, but rather  about the best possible performance of these nonlinear methods of approximation.

Our setting is  a Banach space $X$ equipped with a norm $\|\cdot\|_X$, where we wish to approximate the elements $f$ of 
a compact subset $\cK\subset X$  of $X$ with error measured in this norm. For every fixed $n\in \N$ and every $\gamma\geq 0$, the approximants to $\cK$ will come from  the images $\Phi(y)\in X$ of $\gamma$-Lipschitz maps
$\Phi:(B_{Y_k},\|\cdot\|_{Y_k})\to X$, where $k\leq n$,  and $B_{Y_k}$ is the unit ball in $\R^k$ with respect to some norm $\|\cdot\|_{Y_k}$ in $\R^k$.
 The quality of this approximation is  a critical element in the design and analysis of various numerical methods, among which are DNNs. Note that any numerical method based on  Lipschitz mappings will have performance  no better than the optimal performance of this approximation method.   On the other hand,  it may not be easy to actually design a numerical method for a particular application  that achieves 
 this optimal performance.  
 
  In our analysis, we examine model classes $\cK\subset X$,
 i.e.,  compact subsets $\cK$ of $X$, that  summarize what we know about the target function $f$.   Classical model classes $\cK$ are  finite balls in smoothness spaces like the Lipschitz, Sobolev, or Besov spaces. 
 The {\it Lipschitz} widths  $d^\gamma_n(\cK)_X$ then  quantify the best possible performance of the above approximation methods on a given model class $\cK$.

 The paper is organized as follows. Some of  the basic properties of Lipschitz  widths  are discussed in \S\ref{S2} and \S\ref{S3}, where it is  shown that  for a fixed $n\in \N$, $d_n^\gamma(\cK)_X$ is a  continuous function of 
  $\gamma\geq 0$, see Theorem \ref{continuous}. We also prove the  statements 
 $$
 \lim_{n\to \infty}d_n^\gamma(\cK)_X=0, \quad \hbox{and}\quad \lim_{\gamma\to \infty} d_n^\gamma(\cK)_X=0,
 $$
 each of  which  characterizes the set $\cK$ as a totally bounded set, see Lemma \ref{Lem1} and Lemma \ref{lem0}.
 
The relation between Lipschitz widths and entropy numbers $\e_n(\cK)_X$ is investigated in  \S\ref{S4},   Theorem \ref{TB}, where among other things, 
we show that for any compact subset $\cK\subset X$  of a Banach space $X$ we have
$$
d_n^{2\rad(\cK)}(\cK)_X\leq \e_n(\cK)_X, \quad n=1,2,\dots.
$$
Examples are given to show  that this inequality is almost optimal. We also discuss in this section 
estimates from below and above for the Lipschitz width $d_n^\gamma (\cK)_X$,  provided bounds for the entropy numbers $\e_n(\cK)_X$ are available, see Corollary  \ref{cr1}.
Some of our estimates are optimal,  as demonstrated in Theorem \ref{mainP}, where we show that the Lipschitz widths could be smaller than the entropy numbers for certain compact classes $\cK$.

Since the Lipschitz width is a new concept of width, we compare it with some of the well known classical widths. 
We show that for appropriate values of the parameter $\gamma$,  Lipschitz widths are smaller than the Kolmogorov widths, see \S\ref{S5}, Theorem \ref{TK}. They are also smaller than the stable manifold widths, see \S\ref{S6}, Theorem \ref{TS}. However, as demonstrated by the provided Examples, in both cases, their actual behavior may be very different. 

At last, in \S\ref{S7},  we discuss  the Lipschitz widths $d_n^{\gamma_n} (\cK)_X$  with  $\gamma_n=C'n^\delta\lambda^n$ and show that they provide a theoretical benchmark for the performance of certain DNN approximation, see Theorem \ref{NN}. The  analysis of these  widths is performed   in Theorem \ref{widthsfrombelownew} and Corollary \ref{cNN}, where it is demonstrated that there is indeed a  gain in the performance of the Lipschitz width $d_n^{\gamma_n} (\cK)_X$  when compared to the entropy numbers $\e_n(\cK)_X$ in the following sense
$$
\hbox{if}\quad\e_n(\cK)_X\asymp \frac{[\log_2n]^\beta}{n^\alpha} \quad \Rightarrow\quad 
d_n^{\gamma_n}(\cK)_X\asymp \frac{[\log_2n]^\beta}{n^{2\alpha}}.
$$
This estimate, when applied in the case  of $\cK$ being the unit ball of certain Besov spaces,  extends the results  
 from \cite{DHP} to the case when error is measured in $L_p$, $p\neq \infty$.

\section{Definition and basic properties}
\label{S2}
  We are mainly interested in compact sets, however we define the basic concepts for bounded sets.
We consider  a bounded  subset  $\cK\subset X$  of a Banach space $(X,\|\cdot\|_X)$ with norm $\|\cdot\|_X$ and 
denote by    $(\R^n,\|.\|_{Y_n})$, $n\geq 1$
the $n$-dimensional Banach space with a fixed norm $\|\cdot\|_{Y_n}$. 
For $\gamma\geq 0$, we define 
 the {\it fixed Lipschitz} width
\be
\label{fixed}
d^\gamma(\cK,Y_n)_X:= \inf_{\Phi_n} \sup_{f\in \cK}\inf_{y\in B_{Y_n}}\|f-\Phi_n(y)\|_X,
\ee
 where the infimum is taken over all Lipschitz mappings
 $$
 \Phi_n:(B_{Y_n},\|\cdot\|_{Y_n})\to X,\quad B_{Y_n}:=\{y\in \R^n:\,\,\|y\|_{Y_n}\leq 1\},
 $$ 
that satisfy  
 the Lipschitz condition 
\be
\label{Lip}
 \sup_{y,y'\in B_{Y_n}}
 \frac{\|\Phi_n(y)-\Phi_n(y')\|_X}{\|y-y'\|_{Y_n}}\leq \gamma,
\ee
 with constant $\gamma$. 
 Next, 
we define the
{\it Lipschitz} width
\be
\label{Lwidth}
d_n^\gamma(\cK)_X:=\inf_{k\leq n}\inf_{\|\cdot\|_{Y_k}}d^\gamma(\cK,Y_k)_X,
\ee
where the infimum is taken over  all norms $\|\cdot\|_{Y_k}$ in $\R^k$ and all $k\leq n$.
Clearly, we have that for every norm $\|\cdot\|_{Y_n}$ on $\R^n$,
\be
\label{n111}
 d_n^\gamma(\cK)_X\leq d^\gamma(\cK,Y_n)_X, \quad n\geq 1.
\ee

Before going further, let us recall  the  definition  of a diameter and radius of a bounded set  $\cM\subset X$, 
$$
\diam(\cM):=\sup_{f,g\in \cM}\|f-g\|_X\leq 2\inf_{g\in X}\sup_{f\in \cM}\|f-g\|_X=:2\rad(\cM).
$$
From \eref{Lip} we see that a $0$-Lipschitz function is simply the constant function. Thus,  for any $Y_n$ we get
\be \label{Lip0}
\rad \cK=d_n^0(\cK,Y_n)_X=d_n^0(\cK)_X.
\ee

We next list some elementary properties of the Lipschitz widths $d_n^\gamma(\cK)_X$ that we gather in the following lemma.

 \begin{lemma}
 \label{LL1}
 For any  bounded subset $\cK\subset X$ of a Banach space $X$, any $n\in \N$,  and any $\gamma>0$, we have
 \begin{enumerate}
\item  The Lipschitz width $d_n^\gamma(\cK)_X$ is given by
\be
 \label{Lwidthnew}
 d_n^\gamma(\cK)_X=
\inf_{\|\cdot\|_{Y_n}}d^\gamma(\cK,Y_n)_X.
 \ee
\item{ We can restrict the infimum in \eref{Lwidthnew}  only to normed spaces $(\R^n,\|.\|_{\cY_n})$ with the additional  property that  the norm $\|\cdot\|_{\cY_n}$   satisfyies the 
condition 
	\be
	 \label{Auerbach}
	\|y\|_{\ell_\infty^n}:=\max_j|y_j|\leq \|y\|_{\cY_n}\leq \sum_{j=1}^n|y_j|=:\|y\|_{\ell_1^n}, \quad y=(y_1,y_2,\ldots,y_n)\in \R^n.
	\ee
} 
  \item The space $(\R^n,\|\cdot\|_{Y_n})$ in  \eref{fixed} and \eref{Lwidthnew} can be replaced by any normed space $(X_n,\|\cdot\|_{X_n})$  of dimension $n$, that is
 \be
\label{Lipnew}
d_n^\gamma(\cK)_X=\inf_{\|\cdot\|_{X_n}}d^\gamma(\cK,X_n)_X, \quad \hbox{where} \quad 
d^\gamma(\cK,X_n)_X=\inf_{\Phi_n} \sup_{f\in \cK}\inf_{x\in B_{X_n}}\|f-\Phi_n(x)\|_X,
\ee
with  $B_{X_n}:=\{x\in X_n:\,\,\|x\|_{X_n}\leq 1\}$.
\item  $d_n^{\gamma}(\cK)_X$  is a monotone decreasing function of $\gamma$ and $n$. More precisely,
\begin{itemize}
	\item If $\gamma_1\leq \gamma_2$ then $d_n^{\gamma_2}(\cK)_X\leq d_n^{\gamma_1}(\cK)_X$;
       \item If $n_1\leq n_2$ then $d_{n_2}^{\gamma}(\cK)_X\leq d_{n_1}^{\gamma}(\cK)_X$.
       
       \end{itemize}
       \item  For every fixed $n\in \N$ and $\gamma\geq 0$, we have $d_n^\gamma(\cK)_X\leq \rad(\cK)<\infty$.
       \item For every fixed $n\in \N$ and $\gamma\geq 0$, we have $d_n^\gamma(\cK)_X=d_n^\gamma(\bar \cK)_X$ where $\bar \cK$ denotes the closure of $\cK$.
 \end{enumerate}
 \end{lemma}
 \noindent
 {\bf Proof:} 
Since 
$$
d_n^\gamma(\cK)_X{\leq
		\inf_{\|\cdot\|_{Y_n}}d^\gamma(\cK,Y_n)_X},
		$$
to show (i), it suffices to show that for every  norm $\|\cdot\|_{Y_k}$ with $1\leq k<n$, there exists  a norm $\|\cdot\|_{Y_n}$  on $\R^n$ such that 
\be
\label{wa}
d^\gamma(\cK,Y_n)_X\leq d^\gamma(\cK,Y_k)_X.
\ee
Indeed, let us fix a $\gamma$-Lipschitz map $\Phi_k:(B_{Y_k},\|\cdot\|_{Y_k})\to X$ which achieves  $d^\gamma(\cK,Y_k)_X$ (if such map does not exist, we can use limiting arguments). 
We then define the norm $\|\cdot\|_{Y_n}$ as  
$$
\|(y,y')\|_{Y_n}:=\|y\|_{Y_k}+\|y'\|_{\R^{n-k}},
$$
 where $\|.\|_{\R^{n-k}}$ is any norm in $\R^{n-k}$,  and a mapping $\Phi_n:(B_{Y_n},\|\cdot\|_{Y_n})\to X$ as 
 $$
 \Phi_n((y,y')):=\Phi_k(y).
 $$ 
 Clearly, $\Phi_n$ is a $\gamma$-Lipschitz mapping since
 $$
 \|\Phi_n((y,y'))-\Phi((z,z'))\|_X=\|\Phi_k(y)-\Phi_k(z)\|_X\leq \gamma\|y-z\|_{Y_k}\leq\gamma \|(y,y')-(z,z')\|_{Y_n},
 $$
 and thus \eref{wa} holds.

{
\noindent
Next, we use (i) to prove (ii).  Let $(\R^n,\|\cdot\|_{Y_n})$ be any normed space. It follows from the Auerbach lemma (see e.g. \cite[p.43]{CS} or \cite[II.E.11]{PW}),  that  we can find vectors $(\bar v_j)_{j=1}^n\subset R^n$ and linear functionals $(f_j)_{j=1}^n\subset (\R^n)^*$ such that
	\be
	\label{l1}
	\|\bar v_j\|_{Y_n}=\|f_j\|_{ (\R^n)^*}=1, \quad  j=1,\dots,n,
	\ee
	 and 
	\be
	\label{l2}
	 f_i(\bar v_j)=\delta_{i,j}=\begin{cases}
	 1,\quad i=j,\\
	 0, \quad i\neq j.
	 \end{cases}
	\ee
We  define a new norm $\|.\|_{\cY_n}$ on $\R^n$ as 
$$
\|y\|_{\cY_n}:= \|\sum_{j=1}^n y_j \bar v_j\|_{Y_n},\quad y=(y_1,\ldots,y_n)\in\R^n,
$$
which, using the triangle inequality and \eref{l1}, satisfies the inequality
\be
\label{p1}
\|y\|_{\cY_n}\leq \sum_{j=1}^n |y_j|\|\bar v_j\|_{Y_n}=\sum_{j=1}^n |y_j|.
\ee
On the other hand, using \eref{l1} and \eref{l2}, we have
\be
\label{p2}
\displaystyle{\|y\|_{\cY_n}= \|\sum_{j=1}^n y_j \bar v_j\|_{Y_n}=\sup_{f\in (\R^n)^*, \,\,\|f\|_ {(\R^n)^*}=1}
|f(\sum_{j=1}^n y_j\bar v_j)|\geq |f_i(\sum_{j=1}^n y_j\bar v_j)|=|y_i|, \quad i=1,\ldots,n.}
\ee
Therefore, it follows  from \eref{p1} and \eref{p2} that the newly defined norm satisfies \eref{Auerbach}. If we consider the mapping $\phi_0$ defined as
$$
\phi_0(y):=\sum_{j=1}^ny_j\bar v_j,\quad y=(y_1,\ldots,y_n)\in R^n,
$$
one can show   that $\phi_0:(B_{\cY_n},\|\cdot\|_{\cY_n})\to (B_{Y_n},\|\cdot\|_{Y_n})$ and that $\phi_0(B_{\cY_n})=B_{Y_n}$.
Now, for any $\gamma$-Lipschitz mapping
$\Phi_n:(B_{Y_n},\|\cdot\|_{Y_n})\to X$, we define the map $\tilde \Phi_n:(B_{\cY_n},\|\cdot\|_{\cY_n})\to X$ as
$$
\tilde \Phi_n:=\Phi_n\circ \phi_0.
$$
Note that $\tilde \Phi_n$ is $\gamma$-Lipschitz since
\begin{eqnarray}
\nonumber
\|\tilde \Phi_n(y')-\tilde \Phi_n(y)\|_X&=&\|\Phi_n\circ \phi_0(y')-\Phi_n\circ \phi_0(y)\|_X\leq \gamma\|\phi_0(y')-\phi_0(y)\|_{Y_n}\\
\nonumber
&=&\gamma\|\sum_{j=1}^n(y'_j-y_j)\bar v_j\|_{Y_n}=
\gamma\|y'-y\|_{\cY_n}.
\nonumber
\end{eqnarray}
In addition,  $\tilde\Phi_n(B_{\cY_n})=\Phi_n(B_{Y_n})$, and therefore (ii) follows from (i).
}

 To prove (iii),  we fix a basis $\{\phi_1,\ldots,\phi_n\}\in X_n$, the mapping $\kappa:X_n\to\R^n$ given by
$$
\kappa(g)=(y_1,\ldots,y_n), \quad \hbox{for}\quad g=\sum_{j=1}^ny_j\phi_j,
$$
is isometry between $(X_n,\|\cdot\|_{X_n})$ and  $(\R^n,\|\cdot\|_{Y_n})$, where  $\|x\|_{Y_n}:=\|g\|_{X_n}$. Thus, each norm $\|\cdot\|_{Y_n}$ on $\R^n$ induces a norm 
$\|\cdot\|_{X_n}$ on $X_n$ and vice versa.
Moreover,  the mappings $ \Phi_n\circ \kappa:B_{X_n}\to X$ and $ \Phi_n:B_{Y_n}\to X$ have the same Lipschitz constants, which shows the equivalence of the two definitions
\eref{Lwidthnew}  and \eref{Lipnew}.

Next, (iv)   and (vi)  follow  directly from the definition, and (v) follows from \eref{Lip0} and (iv).
\hfill $\Box$

\subsection{Packing, covering and entropy numbers} 

Before going further, we recall in this section the well known concepts of packing, covering, and entropy numbers for  compact sets $\cM$, which we will use in our study of 
Lipschitz widths. The reader may find a more detailed exposition of those concepts in many books, see, for example,  \cite{CS, AP,LGM}.

{\bf Minimal  $\e$-covering number $ N_\e(\cM)$ of a compact set $\cM\subset X$: } 

\noindent
A collection 
$\{g_1,\ldots,g_{m}\}\subset X$
of  elements of $X$
is called an $\e$-covering  of $\cM$ if 
$$
\cM\subset \bigcup_{j=1}^m B(g_j,\e), \quad \hbox{where}\quad B(g_j,\e):=\{f\in X:\,\|f-g_j\|_X \leq \e \}.
$$
An  $\e$-covering of $\cM$ whose cardinality   is minimal is called {\it minimal $\e$-covering} of $\cM$. 
We denote by  $N_\e(\cM)$   the  cardinality of the  minimal $\e$-covering  of $\cM$.

{\bf  Minimal   inner $\e$-covering number $\tilde N_\e(\cM)$ of a compact set $\cM\subset X$: }

\noindent
 It is defined exactly as $N_\e(\cM)$ 
but we additionally require that the centers  $\{g_1,\ldots,g_{m}\}$ of the covering are elements from $\cM$.

{\bf Entropy numbers $\e_n(\cM) _X$ of a compact set $\cM\subset X$: } 

\noindent
For every fixed $n\geq 0$, the {\it entropy number} $\e_n(\cM) _X$ is the infimum of all $\e>0$ for which $2^n$ balls with centers from $X$ and radius $\e$ cover $\cM$.  
If we put the additional restriction that the centers of these balls are from $\cM$, then we define the so called {\it   inner entropy number} $\tilde \e_n(\cM)_X$.
Formally,  we  write
$$ 
\e_n(\cM)_X=\inf\{ \e>0 \ :\ \cM \subset \bigcup_{j=1}^{2^n} B(g_j,\e), \ g_j\in X, \ j=1,\ldots,2^n\},
$$
$$
\tilde \e_n(\cM)_X=\inf\{ \e>0 \ :\ \cM \subset \bigcup_{j=1}^{2^n} B(h_j,\e), \ h_j\in \cM, \ j=1,\ldots,2^n\}.
$$

{\bf Maximal $\e$-packing  number $\tilde P_\e(\cM)$ of a compact set $\cM\subset X$: } 

\noindent
A collection 
$\{f_1,\ldots,f_{\ell}\}\subset \cM$
of  elements from $\cM$
is called an $\e$-packing of $\cM$ if 
$$
\min_{i\neq j}\|f_i-f_j\|_X>\e.
$$
An  $\e$-packing of $\cM$ whose size   is maximal is called {\it maximal $\e$-packing} of $\cM$.  We denote by $\tilde P_\e(\cM)$ the cardinality of the  maximal $\e$-packing of $\cM$.

We have the following inequalities for every $\e>0$ and every compact set $\cM$
\be 
\label{known}
\tilde P_\e(\cM)\geq \tilde N_\e(\cM)\geq \tilde P_{2\e}(\cM),
 \ee
\be \label{innerentropy}
\e_n(\cM)_X\leq \tilde \e_n(\cM)_X\leq 2 \e_n(\cM)_X.
\ee

\begin{remark}
\label{totally bounded}  
Let us recall the classical relations between those concepts and compactness.
We call the set $\cM$ {\it totally bounded } if for every $\e>0$ we have $N_\e(\cM)<\infty$. This is equivalent to the fact that  $\lim_{n\to \infty}\e_n(\cM)_X=0$. Each compact set is totally bounded. Actually a subset $\cM$ of a Banach space is compact if and only if it is totally bounded and closed. The interested reader will find a detailed study  on the topic in many books on functional analysis or metric topology.
\end{remark}

\begin{remark}
\label{zeroremark}
In what follows later, we will use the fact that the Lipschitz widths and the entropy numbers are invariant with respect to translation, that is,  for any $n\in \N, \gamma\geq 0$, and any $f\in X$ we have 
$$
d_n^\gamma(\cK)_X=d_n^\gamma(\cK-f)_X, \quad 
\e_n(\cK)_X=\e_n(\cK-f)_X.
$$
\end{remark}

\subsection{Dependence of $d_n^\gamma(\cK)_X$ on $\gamma$}

We start this section by proving the fact  that the Lipschitz width $d_n^\gamma(\cK)_X$ is a continuous function of $\gamma$. To do that, we first prove the 
following lemma.
\begin{lemma}
\label{LP}
 For  every $n\geq 1$, every $\gamma> 0$,  and  every norm 
$\|\cdot\|_{Y_n}$ in $\R^n$, the fixed Lipschitz width $d^\gamma(\cK,Y_n)_X$ satisfies the inequality
\be
\label{eq1}
\rad(\cK)-\gamma\leq d^\gamma(\cK,Y_n)_X\leq \rad(\cK).
\ee
\end{lemma}
\noindent
{\bf Proof:}
If we fix $g\in X$ and take $\Phi(y)=g$ for every $y\in B_{Y_n}$, we have that $\Phi$ is $\gamma$-Lipschitz for every $\gamma>0$ and thus
$$
d^\gamma(\cK,Y_n)_X\leq \sup_{f\in\cK}\|f-g\|_X, 
$$
which gives $d^\gamma(\cK,Y_n)_X\leq\rad(\cK)$.
 To show the left hand-side inequality in   \eref{eq1}, we notice that for any $\gamma$-Lipschitz map $\Phi$, every $f\in \cK$ and $y\in B_{Y_n}$  we have 
$$
\|f-\Phi(y)\|_X\geq \|f-\Phi(0)\|_X-\|\Phi(0)-\Phi(y)\|_X\geq \|f-\Phi(0)\|_X-\gamma,
$$ 
since $\|\Phi(0)-\Phi(y)\|_X\leq \gamma\|y\|_{Y_n}\leq \gamma$.
Therefore we obtain the inequality
$$
\sup_{f\in \cK} \inf_{y\in B_{Y_n}}\|f-\Phi(y)\|\geq \sup_{f\in \cK} \|f-\Phi(0)\|-\gamma,
$$
which  gives
\be
\label{q1}
d^\gamma(\cK,Y_n)_X\geq  \inf_\Phi \sup_{f\in \cK}\|f-\Phi(0)\|_X-\gamma.
\ee
Note now that for every $\Phi$, 
$$
\sup_{f\in \cK}\|f-\Phi(0)\|_X\geq \inf_{g\in X}\sup_{f\in \cK}\|f-g\|_X=\rad(\cK),
$$
and thus it follows from \eref{q1} that
$$
d^\gamma(\cK,Y_n)_X\geq \rad(\cK)-\gamma,
$$
and the proof is completed.
\hfill $\Box$

\begin{theorem}\label{continuous}
For every compact subset $\cK\subset X$  of a Banach space $X$ and any $n\in \N$,  
the Lipschitz width $d_n^\gamma(\cK)_X$ is a continuous function of $\gamma  \geq 0$.	
\end{theorem}
\noindent
{\bf Proof:} 
We first show the continuity of the Lipschitz width at  $\gamma=0$.
It follows from  Lemma \ref{LP} and Lemma \ref{LL1}, (i) that
$$
\rad(\cK)-\gamma\leq d_n^\gamma(\cK)_X\leq \rad(\cK).
$$
 We let $\gamma\rightarrow 0$ 
and obtain
$$
\lim_{\gamma\to 0} d_n^\gamma(\cK)_X=\rad(\cK)=d_n^0(\cK)_X,
$$
which proves the continuity at $\gamma=0$, see \eref{Lip0}.

To show that the Lipschitz width is continuous for $\gamma>0$, we fix $n\in \N$ and
 denote by 
$$
h(\gamma):=d_n^\gamma(\cK)_X.
$$
According to Lemma \ref{LL1},  (v), 
$h(\gamma)<\infty$ for every $\gamma>0$.  
Let us  assume that $h$ is not a continuous function. Then, there exist   $\gamma_0>0$, $\delta>0$,  and a sequence of positive numbers $\e_k\to 0$, such that 
\be \label{continuous0}
 h(\gamma_0+\e_k)+\delta\leq h(\gamma_0-\e_k), \quad \hbox{for every}\,\, k.
\ee
We fix $\e:=\e_k<\gamma_0$. From the definition of Lipschitz widths, there exists a $(\gamma_0+\e)$-Lipschitz map $\Phi_n:(B_{Y_n}, \|\cdot\|_{Y_n})\rightarrow X$ such that
\be \label{continuous1} 
h(\gamma_0+\e) \leq \sup_{f\in \cK} \inf_{y\in B_{Y_n}} \|f-\Phi_n(y)\|_X\leq h(\gamma_0+\e)+\e.
\ee
Now we define the mapping 
$$
\tilde \Phi_n:= \xi\Phi_n, \quad \hbox{where}\quad 
\xi:=\frac{\gamma_0-\e}{\gamma_0+\e}\quad \hbox{and}\quad  0<\xi<1.
$$
Clearly, $\tilde \Phi_n$ is  a $(\gamma_0-\e)$-Lipschitz mapping, and therefore
\begin{eqnarray*}
h(\gamma_0-\e)&\leq &\sup_{f\in \cK}\inf_{y\in B_{Y_n}}\|f-\tilde \Phi_n(y)\|_X= \sup_{f\in \cK}\inf_{y\in B_{Y_n}}	\|\xi(f-\Phi_n(y))+(1-\xi)f\|_X\\
&\leq& \sup_{f\in \cK}\inf_{y\in B_{Y_n}}\left(\xi\|f-\Phi(y)\|_X+(1-\xi)\|f\|_X\right) \\
&\leq& \xi \sup_{f\in \cK}\inf_{y\in B_{Y_n}} \|f-\Phi_n(y)\|_X +(1-\xi)\sup_{f\in \cK}\|f\|_X\\
 &\leq& \xi ( h(\gamma_0+\e)+\e) +(1-\xi)C,
\end{eqnarray*}	
where we have used \eref{continuous1} and the fact that $\sup_{f\in \cK}\|f\|_X=C<\infty $ (since $\cK$ is compact). The latter inequality 
 and \eref{continuous0} give 
$$
 h(\gamma_0+\e)+\delta\leq   h(\gamma_0-\e) \leq \xi (h(\gamma_0+\e)+\e) +(1-\xi)C
 $$
which is a contradiction for a sufficiently small  $\e=\e_k$ since $\xi\to 1$ as $\e_k\to 0$.
\hfill $\Box$

We finish the investigation of the behavior of the Lipschitz width with respect to $\gamma$ with the following lemma.

\begin{lemma}
\label{Lem1}
For any $\cK\subset X$,  the set  $\cK$ is 
 totally bounded  
iff for every $n\geq 1$ 
		$$
		\lim_{\gamma \to \infty} d_n^\gamma(\cK)_X=0.
		$$ 
\end{lemma}
\noindent
{\bf Proof:} Assume that $\cK$ is  totally bounded. From the monotonicity of the Lipschitz width with respect to $n$, see Lemma \ref{LL1}, (iii), it suffices to consider only the case $n=1$. For $\delta >0$, we  fix a  minimal  delta covering $(f_j)_{j=1}^{\cN_\delta(\cK)}$  of $\cK$ and choose 
$\gamma $ such that 
 $$
 2\gamma{ \geq}\diam \cK\cdot (\cN_\delta(\cK)-1). 
 $$
We consider the  points 
$$
 t_j:=-1+ 2\frac{j-1}{N_\delta(\cK)-1}, \quad j=1,\dots, \cN_\delta(\cK),
 $$
in the unit ball of  $(\R,|.|)$,  that is  $([-1,1],|\cdot|)$,  and  define $\Phi:[-1,1]\rightarrow X$ as the continuous piecewise linear function such that 
  $$
  \Phi(t_j)=f_j, \quad j=1,\dots, \cN_\delta(\cK).
  $$
Its Lipschitz constant
is no more than
$$
 \max_{j=1,\dots,\cN_\delta(\cK)-1}\frac{\|f_{j+1}-f_{j}\|_X}{|t_{j+1}-t_j|}\leq \frac{\diam \cK\cdot (\cN_\delta(\cK)-1)}{2}\leq \gamma.
$$
and we have 
$$
   \sup_{f\in\cK}\inf_{y\in [-1,1]}\|f-\Phi(y)\|_X\leq \delta.
 $$
This gives 
$$
d_1^\gamma(\cK)_X\leq \delta,
$$
and therefore $\lim_{\gamma\to\infty}d_1^\gamma(\cK)_X=0$.

We now  fix $n\geq 1$ and  show that $\lim_{\gamma\to\infty}d_n^\gamma(\cK)_X=0$ implies that $\cK$ is  totally bounded.
We prove it  by showing  that if $\cK$ is not  totally bounded, we can find  $\delta>0$ such that 
$d_n^\gamma(\cK)_X\geq \delta$ 
for every $\gamma>0$.
So, we now  assume that $\cK$ is not  totally bounded, which
 implies that there exists $\delta_0>0$ and an infinite $\delta_0$-packing set that we will denote $(h_j)_{j=1}^\infty$.

 Let us fix a norm $\|\cdot\|_{Y_n}$ on $\R^n$ and consider any $\gamma$-Lipschitz map
$\Phi:(B_{Y_n},\|\cdot\|_{Y_n})\to X$. 
 We then take $\epsilon <\delta_0/ 3\gamma$ and denote by 
 $\{y_j\}_{j=1}^N\subset B_{Y_n}$
 a finite  $\epsilon$-covering  of $B_{Y_n}$.

Note that at most one $h_i$ can belong to any of the sets $B(\Phi(y_j),\delta_0/ 2)$. Indeed,   if  we assume  that $h_{j_1}\neq h_{j_2}$ and $h_{j_1},h_{j_2}\in B(\Phi(y_j),\delta_0/ 2)$ for some $j\in \{1,\ldots,N\}$, then 
 $$
\|h_{j_1}-h_{j_2}\|_X\leq \|h_{j_1}-\Phi(y_j)\|_X+\|\Phi(y_j)-h_{j_2}\|_X \leq \delta_0,
 $$
which  contradicts the fact that  $\{h_j\}_{j=1}^\infty$ is a $\delta_0$-packing set.
Therefore, there exists $s\in \{1,2,\ldots\}$ such that 
\be
\label{nn}
h_s\notin \bigcup_{j=1}^N B(\Phi(y_j),\delta_0/ 2).
\ee
We also know that for every $y\in B_{Y_n}$ there is $j^*\in \{1,2,,\ldots,N\}$ such that 
 $\|y-y_{j^*}\|_{Y_n} \leq\e$, and therefore
 \be
\label{lq}
\|\Phi(y)-\Phi(y_{j^*})\|_X\leq \gamma\|y-y_{j^*}\|_{Y_n} \leq \gamma\e<\delta_0/3.
\ee
Using \eref{lq}, \eref{nn}  and the triangle inequality, we obtain
$$
\|h_s-\Phi(y)\|_X\geq \|h_s-\Phi(y_{j^*})\|_X-\|\Phi(y)-\Phi(y_{j^*})\|_X\geq \delta_0/2-\delta_0/3=\delta_0/6,
$$
which gives
$$
d^\gamma(\cK,\|\cdot\|_{Y_n})_X\geq \delta_0/6, 
$$
and therefore
\be
\label{noncompact}
 d_n^\gamma(\cK)_X\geq \delta_0/6.
\ee
Notice that the choice of  $\gamma$ is arbitrary, so we have the above inequality for every $\gamma$, and the proof is completed. 
\hfill$\Box$

\begin{remark}
Note that all statements  in this paper are valid for sets $\cK$ whose closures are compact rather than sets $\cK$ that are compact. Therefore, since we work in Banach spaces,  all statements are valid 
for 
$\cK$ being  only a totally bounded set rather than a compact set.
\end{remark}

\begin{remark}
It follows from the proof of Lemma {\rm \ref{Lem1}},  see \eref{noncompact},  that   if $\cK$ is not totally bounded  then there exists $\delta>0 $ such that
	$$
	 \inf _{\gamma>0,n>0}d_n^\gamma(\cK)_X\geq \delta.
	 $$
\end{remark}

\subsection{Dependence of the Lipschitz width $d_n^\gamma(\cK)_X$ on $n$.}

In this section  we discuss the behavior of the Lipschitz width with respect to $n$.
 The following Lemma holds.
\begin{lemma} 
\label{lem0}
Let $\cK\subset X$ be a subset of a Banach space $X$. If there exists $\gamma>0$ such that  
$\lim_{n\to \infty} d_n^\gamma (\cK)_X=0$, then $\cK$ is totally bounded (i.e. its closure is compact).
\end{lemma}
\noindent
{\bf Proof:} To prove the lemma, we   fix $\eta>0$ and  show that $\cK$ is contained in the union of a finite 
 collection of balls with radius $\eta$.
Since $\lim_{n\to \infty} d_n^\gamma (\cK)_X=0$, we can find
an integer $n_0$ such that 
$$
d_{n_0}^\gamma(\cK)_X<\eta/2,
$$
and therefore there exists a 
norm $\|\cdot\|_{Y_{n_0}}$ in $\R^{n_0}$ and a $\gamma$-Lipschitz map $\Phi:(B_{Y_{n_0}},\|\cdot\|_{Y_{n_0}})\rightarrow X$ such that
$$
	\sup_{f\in \cK}\inf_{y\in B_{Y_{n_0}}}\|f-\Phi(y)\|<\eta/2.
$$
More precisely, for every $f\in \cK$, we can find $y\in B_{Y_{n_0}}$ such that
\begin{equation}
\label{J12;1}
\|f-\Phi(y)\|_X<\eta/2.
\end{equation}
 Let $\{y_j\}_{j=1}^N\subset B_{Y_{n_0}}$ be an $\eta/(2\gamma)$-covering for the compact set $B_{Y_{n_0}}$,
that is
$$
B_{Y_{n_0}}\subset \bigcup_{n=1}^N B(y_j, \eta/(2\gamma)),
$$
and therefore  for every $y\in B_{Y_{n_0}}$ we can find $y_j $, $j\in \{1,2,\ldots,N\}$, such
 that $\|y-y_j\|_{Y_{n_0}} \leq \eta/(2\gamma)$. Thus we have
$$
\|\Phi(y)-\Phi(y_j)\|_X\leq \gamma\|y-y_j\|_{Y_{n_0}}\leq\eta/2,
$$
and
$$
\Phi(B_{Y_{n_0}})\subset \bigcup_{n=1}^N B(\Phi(y_j), \eta/2).
$$
From the    latter result  and \eref{J12;1} it follows that 
$\cK\subset \bigcup_{n=1}^N B(\phi(y_j), \eta)$, and the proof is completed.
 \hfill $\Box$

The converse statement of Lemma \ref{lem0} is also true, see Corollary \ref{tend0} and Corollary \ref{tend0_1}.

 \section{Further properties of Lipschitz widths}
\label{S3}

\subsection{Properties of Lipschitz mappings}

Before focusing our attention on the Lipschitz widths, we want to state and prove a lemma which shows the behavior of the entropy numbers of an image of a  $\gamma$-Lipschitz mapping. More precisely, 
the following holds. 

\begin{lemma} 
\label{DD}
Consider the two normed spaces $(X_0,\|\cdot\|_{X_0})$ and  $(X_1,\|\cdot\|_{X_1})$ and the $\gamma$-Lipschitz map
$$
\Phi:(\cK_0,\|\cdot\|_{X_0})\to (\cK_1,\|\cdot\|_{X_1}), \quad \cK_0\subset X_0, \quad \cK_1\subset X_1.
$$
Then the following holds:
\begin{enumerate}
\item if $\Phi(\cK_0)=\cK_1$, then  
\be
\label{EE}
  \tilde\e_k(\cK_1)_{X_1}\leq \gamma \tilde\e_k(\cK_0)_{X_0}, \quad k=1,2,\ldots.
\ee
In particular, if $\Phi:(B_{Y_n},\|\cdot\|_{Y_n})\to
(B_{Z_m},\|\cdot\|_{Z_m})$ is a $\gamma$-Lipschitz map from the unit ball  $B_{Y_n}$ onto the unit ball $B_{Z_m}$, then $n\geq m$.
\item if $\Phi(\cK_0)$ approximates $\cK_1$ with accuracy  $\e_2$ and $A\subset \cK_0$ approximates $\cK_0$ with accuracy $\e_1$, then
			$\Phi(A)$ approximates $\cK_1$ with accuracy $ \gamma \e_1 +\e_2$.

\end{enumerate}
\end{lemma}
\noindent
{\bf Proof:} We first prove (i). We consider  the set $\{g_j\}_{j=1}^{2^k}\subset \cK_0$  of $2^k$ elements of $\cK_0$,
 $k\geq 1$,  that are the centers of $2^k$ balls of radius 
$\e\geq \tilde\epsilon_k(\cK_0)_{X_0}$ that cover $\cK_0$. Let 
$$
f_j:=\Phi(g_j), \quad 
j=1,\ldots,2^k,
$$
be the images of these $g_j$'s under $\Phi$.
Then, since $\Phi(\cK_0)=\cK_1$, for every $f\in \cK_1$ there is $g\in \cK_0$ such that $\Phi(g)=f$, and index $j^*\in\{1,\ldots,2^k\}$ such that $ \|g-g_{j^*}\|_{X_0} \leq \e$. Therefore we have
$$
\|f-f_{j^*}\|_{X_1}=\|\Phi(g)-\Phi(g_{j^*})\|_{X_1}\leq \gamma \|g-g_{j^*}\|_{X_0}\leq\gamma \e,
$$
which shows that $\{f_j\}_{j=1}^{2^k}$ provides a covering for $\cK_1$ with radius $ \leq\gamma\e$, and thus
$$
\tilde\e_k(\cK_1)_{X_1}\leq \gamma \e.
$$
 The latter inequality is true for any  $\e$ being the radius of a set of $2^k$ balls that cover  $ \cK_0$, and therefore
 \eref{EE} holds by taking the infimum over all such $\e$.

In the case when $\cK_0=B_{Y_n}$ and $\cK_1=B_{Z_m}$, we derive from  \eref{EE}  and \eref{innerentropy} that  
$$
\e_k(B_{Y_n})\geq (2\gamma)^{-1}\e_k(B_{Z_m}), \quad k=1,2,\dots.
$$
We know  from  (1.1.10) in \cite{CS} 
 that for 
any unit ball $B_{X_\ell}$ of any Banach space $(X_\ell, \|\cdot\|_{X_\ell})$ of dimension $\ell$  we have
$4\cdot 2^{-k/\ell}\geq \e_k(B)\geq 2^{-k/\ell}$,  $k=1,2,\dots $. Thus we get
$$ 
4 \cdot 2^{-k/n}\geq \e_k(B_{Y_n})\geq  (2\gamma)^{-1}\e_k(B_{Z_m})\geq  (2\gamma)^{-1}2^{-k/m}, \quad k=1,2,\dots,
$$
which  can hold only when $n\geq m$.

To show (ii), we take  $f\in \cK_1$, the corresponding  $g\in \cK_0$ such that $\|\Phi(g)-f\|_{X_1}\leq \e_2$ and $g_0\in A$ such that $\|g_0-g\|_{X_0}\leq \e_1$. So, we have
		$$
		\|\Phi(g_0)-f\|_{X_1}\leq \|\Phi(g_0)-\Phi(g)\|_{X_1}+\|\Phi(g)-f\|_{X_1}\leq \gamma\e_1+\e_2,
		$$
		and the proof is completed.

\hfill $\Box$

\subsection{A single norm defines the Lipschitz width}

In this section, see Theorem \ref{specialnorm},  we extend Lemma \ref{LL1},  (ii) and prove that in the definition of Lipschitz width the infimum over all norms is achieved for some norm that satisfies 
\eref{Auerbach}. We use the following
 version of Ascoli's theorem, whose proof can be found in \cite{CDPW},
and which we state below.
\begin{lemma}
	\label{Ascolii}  Let 
	$(X,d)$ be a separable metric space and $(Y,\rho)$ be  a metric space
	for which  every closed ball is compact. Let $F_j:X\to Y$ be a
	sequence of $\gamma$-Lipschitz maps  for which  there exists $a\in X$
	and $b\in Y$ such that $F_j(a)=b$ for $j=1,2,\dots$. Then, there exists
	a subsequence $F_{j_k}$, $k\ge 1$,  which is  point-wise convergent to a function $F:X\to Y$ and
	$F$ is $\gamma$-Lipschitz.  If $(X,d)$  is also compact, then the convergence is uniform.
\end{lemma}
Now we are ready to state and prove the following fact.
\begin{theorem}
\label{specialnorm}
For any   $n\in \N$, any compact set $\cK\subset X$,  and any constant $\gamma>0$ there is a 
norm $\|\cdot\|_Y$ on $\R^n$ satisfying {\rm \eref{Auerbach}}  such that 
$$
d_n^\gamma(\cK)_X=d_n^\gamma(\cK,Y)_X.
$$
\end{theorem}
\noindent
{\bf Proof:} It follows from Lemma \ref{LL1}, (i), (ii) that we can find a sequence $(\Psi^j)_{j=1}^\infty$ of $\gamma$-Lipschitz maps
$\Psi^j:(B_{\cY_j}, \|\cdot\|_{\cY_j})\rightarrow X$, where the norms 
$\|.\|_{\cY_j} $ on $\R^n$ satisfy \eref{Auerbach},  such that 
	$$
	d_j:=\sup_{f\in \cK}\inf_{y\in B_{\cY_j}} \|f-\Psi^j(y)\|_X\rightarrow d_n^\gamma(\cK)\quad \hbox{as}\quad j\rightarrow \infty.
	$$	
There is  a subsequence  $\|.\|_{\cY_{j_k}}$ of the sequence of norms $\|.\|_{\cY_{j}}$ that converges point-wise on $\R^n$ and uniformly on 
 $B_{\ell_\infty^n}$ to a
norm $\|.\|_Y$ on $\R^n$ satisfying \eref{Auerbach}.
Indeed, one can check that the functions $F_j:(\R^n,\|\cdot\|_{\ell_1^n})\to \R$,  defined as $F_j(y):=\|y\|_{\cY_j}$ satisfy $F_j(0)=0$, and
$$
 |F_j(y')-F_j(y)|=|\|y'\|_{\cY_j}-\|y\|_{\cY_j}|\leq \|y'-y\|_{\cY_j}\leq \|y'-y\|_{\ell_1^n},
$$
where we have used \eref{Auerbach}. Thus, the sequence $(F_j)_{j=1}^\infty$ satisfies the conditions of  Lemma \ref{Ascolii} with $\gamma=1$, $a=b=0$, and so
we can find a subsequence $F_{j_k}$ that converges point-wise on $\R^n$ and uniformly on  $B_{\ell_\infty^n}$.
In fact,  the limit function $F$  of this subsequence is 
a norm, which we denote by $\|\cdot\|_Y$. Clearly, this norm satisfies inequalities \eref{Auerbach}.

Now,
 passing to a subsequence,  we will assume that  $\|.\|_{\cY_j}$ converge uniformly on $B_{\ell_\infty^n}$ to the function $\|.\|_Y$.  Thus,   there is $j_0\in \N$ such that for any 
   $j\geq j_0$ there is $\e_j$ with the properties   $0<\e_j<1$ , $\lim_{j\to \infty}\e_j=0$ and 
$$
\| y\|_{\cY_j}
-\e_j \leq \|y\|_Y\leq \|y\|_{\cY_j}+\e_j, \quad \hbox{for all}\quad \|y\|_{\ell_\infty^n}\leq 1.
$$
For example, we can  take 
$$
\e_j:=\sup_{y:\|y\|_{\ell_\infty^n}\leq 1}|\|y\|_{\cY_j}-\|y\|_Y|,
$$
 and $j_0$ big enough.
Since $B_{\cY_j}\subset B_{\ell_\infty^n}$, $j=1,2,\ldots$, and  $B_Y\subset B_{\ell_\infty^n}$,
we have  for all   $y\in B_{\cY_j}\cup B_{Y}$ 
\be
\label{n1}
\| y\|_{\cY_{j}}
-\e_j \leq \|y\|_Y\leq \|y\|_{\cY_{j}}+\e_j.
\ee
 The latter inequality gives that for  $y\in B_Y$ we have
$$
\| y\|_{\cY_{j}}
 \leq 1+\e_j\quad \Rightarrow\quad y\in (1+\e_j)B_{\cY_{j}},
 $$
 and so
\be
\label{incl1}
  B_Y\subset (1+\e_j)B_{\cY_{j}}.
\ee
Next,  let $j\geq j_0$. For any  $y$ with   $\|y\|_{\cY_{j}}\leq 1-\e_j<1$,  we have from \eref{n1} that 
$$
\| y\|_{Y}
 \leq 1\quad \Rightarrow\quad y\in B_{Y},
 $$
 and therefore
\be
\label{incl2}
 (1-\e_j)^{-1}B_{\cY_{j}}\subset B_Y.
\ee
It follows from \eref{incl1} and \eref{incl2}  that 
 \be
\label{ff}
(1-\e_{j})^{-1}B_{\cY_{j}}\subset B_{Y}\subset (1+\e_{j}) B_{\cY_{j}},\quad  j\geq j_0.
\ee
Let us now define the mapping $\tilde\Psi^{j}:(1+\e_j)B_{\cY_{j}}\to X$,    as
$$
\tilde\Psi^{j}(y):=\Psi^{j}((1+\e_j)^{-1}y).
$$
Note that
$$
\|\tilde \Psi^{j}(y')-\tilde \Psi^{j}(y)\|_X    \leq \frac{\gamma}{1+\e_j}\|y'-y\|_{\cY_{j}}<\gamma\|y'-y\|_{\cY_{j}}, \quad y',y\in (1+\e_j)B_{\cY_{j}},
$$
where we have used that $\Psi^j$ is  $\gamma$-Lipschitz. We  denote by $\bar \Psi^j$ the restriction of $\tilde \Psi^j$ on $B_Y$, see \eref{ff}.

Now we fix $f\in \cK$   and  $j\geq j_0$.  For every $\e>0$, we can find $y=y(f,j,\e)\in B_{\cY_j}$ such that $\|f-\Psi^j(y)\|_X<d_j+\e$. We set
$$
z:=y/(1-\e_j)\in (1-\e_j)^{-1}B_{\cY_j}\subset B_Y,
$$
and observe that 
 \begin{eqnarray*}
\inf_{x\in B_Y}\|f-\bar \Psi^j(x)\|_X&\leq& \|f-\bar \Psi^j(z)\|_X=\|f-\tilde \Psi^j(z)\|_X=\|f-\Psi^j((1+\e_j)^{-1}z)\|_X \\
&=&\|f-\Psi^j((1-\e_j^{2})^{-1}y)\|_X 
\leq \|f-\Psi^j(y)\|_X+\|\Psi^j(y)-\Psi^j ((1-\e_j^{2})^{-1}y)\|_X\\
&<& d_j+\e+\gamma\frac{\e_j^2}{1-\e_j^2} \|y\|_{\cY_j}	\leq  d_j+\e+\gamma\frac{\e_j^2}{1-\e_j^2}.
 	\end{eqnarray*}
By letting $\e\to 0$ and taking supremum  over $f\in\cK$, we obtain
$$
d^\gamma(\cK,Y)_X\leq\sup_{f\in \cK}\inf_{x\in B_Y}\|f- \bar \Psi^j(x)\|_X\leq  d_j+\gamma\frac{\e_j^2}{1-\e_j^2}.
$$
Since $d_j\to d_n^\gamma(\cK)_X$ and $\e_j \to 0$ as  $j\to \infty$, we derive that
$d^\gamma(\cK,Y)_X\leq d_n^\gamma(\cK)_X$, and the proof is completed.
\hfill $\Box$

\section{Lipschitz widths and entropy numbers}
\label{S4}

In this section we discuss the relation between the Lipschitz widths $d_n^\gamma(\cK)_X$ and the entropy numbers $\e_n(\cK)_X$ of a compact set $\cK$.

\subsection{ Lipschitz widths are smaller than entropy numbers}
We start with the construction of a particular Lipschitz function that can be viewed as  a sum of `bumps', each one supported on a closed ball from a Banach space $Y$. We use this function 
to show that the Lipschitz widths of a compact set $\cK\subset X$ are smaller than the entropy numbers of that set. {
\begin{lemma}	\label{map}
Let  $(B^j):=(B(y_j,\rho_j))$ be a family of  disjoint open balls 
 in a Banach space $Y$. Then the following holds:
 
{\rm (i)}   for every sequence $(\varphi_j)$ of $\gamma_j$-Lipschitz mappings 
$\varphi_j:Y \to X$, $ j=1,2,\ldots$,  with the property   that $\varphi_j\equiv 0$ on the complement of $B^j$,
the mapping $\Phi:Y \to X$, defined as
\be
\label{hk}
  \Phi_0=\sum_{j} \varphi_j
\ee
 is a Lipschitz map with Lipschitz constant $sup_j\gamma_j$. 
 
 {\rm (ii)}  for any sequence $(f_j)$ of elements $f_j\in X$ with  $\|f_j\|_X= 1$ and any sequence
 $(\sigma_j)$  of real numbers, the mappings $\phi_j:Y\to X$, $j=1,2,\ldots$,  defined as 
$$
\phi_j(y)=\sigma_j\left( 1-\frac{\|y_j-y\|_Y}{\rho_j}\right)_+\cdot f_j, \quad \hbox{where}\quad (t)_+:=\max\{0,t\}, \,\, t\in \R,
$$
are $|\sigma_j|/\rho_j$-Lipschitz mappings. Their sum, the mapping 
$$
  \Phi:=\sum_{j} \phi_j
  $$ 
is a Lipschitz mapping with  Lipschitz constant $\displaystyle{\sup_j |\sigma_j|/\rho_j}$ and  $\Phi(y_j)=\sigma_j f_j$, $j=1,2,\ldots$.
\end{lemma}
{\bf Proof:}  To show (i), we denote by $\gamma:=\sup_j\gamma_j$ and consider several cases:
\begin{itemize}
\item if  $y,y'\in B^j$, then 
$$
\|\Phi_0(y)-\Phi_0(y')\|_X=\|\varphi_j(y)-\varphi_j(y')\|_X\leq \gamma_j\|y-y'\|_Y\leq \gamma \|y-y'\|_Y.$$
\item if both $y,y'$ are outside each of the balls $B^j$, we have 
$\|\Phi_0(y)-\Phi_0(y')\|_X=0$.
\item if $y\in B^j$ and $y'$ is outside the union $\bigcup_kB^k$, we denote by $y_0$ the intersection of $\partial B^j$ and the line segment connecting $y$ with $y'$.
In this case $y_0=sy+(1-s)y'$ for some $s\in [0,1]$.
Then, we have
$\Phi_0(y')=\Phi_0(y_0)=0=\varphi_j(y_0)$, $\Phi_0(y)=\varphi_j(y)$, and 
$$
\|\Phi_0(y)-\Phi_0(y')\|_X=\|\varphi_j(y)-\varphi_j(y_0)\|_X\leq \gamma_j \|y-y_0\|_Y=\gamma_j (1-s)\|y-y'\|_Y\leq \gamma\|y-y'\|_Y.
$$
 \item if $y\in B^j$ and $y'\in B^k$, $j\neq k$, we denote by $y_0$ the intersection of $\partial B^j$ and the line segment connecting $y$ with $y'$, and by $y_1$ the intersection of $\partial B^k$ and the line segment connecting $y$ with $y'$.
 Then, we have 
\begin{eqnarray}
\nonumber
 y_0=s_0y+(1-s_0)y_1, \quad s_0\in [0,1], 
 \nonumber \\
 y_1=s_1y+(1-s_1)y', \quad s_1\in [0,1], 
 \nonumber
\end{eqnarray}
Moreover,  $\Phi_0(y)=\varphi_j(y)$, $\Phi_0(y')=\varphi_k(y')$,  $\varphi_j(y_0)=0=\varphi_k(y_1)$, and therefore
\begin{eqnarray}
\nonumber
\|\Phi_0(y)-\Phi_0(y')\|_X&=&\|\varphi_j(y)-\varphi_k(y')\|_X\leq \|\varphi_j(y)-\varphi_j(y_0)\|_X+\|\varphi_k(y_1)-\varphi_k(y')\|_X 
\end{eqnarray}
For each pair $y$, $y_0$ and $y_1$, $y^\prime$  we apply the previous case and get
\begin{eqnarray*}
\|\Phi_0(y)-\Phi_0(y')\|_X &\leq& \gamma(\|y-y_0\|_Y+\|y_1-y^\prime \|_Y)
 =\gamma((1-s_0)\|y-y_1\|_Y+s_1\|y-y'\|_Y)\\
&=&\gamma((1-s_0)(1-s_1)\|y-y'\|_Y+s_1\|y-y'\|_Y)\\
 &\leq& \gamma((1-s_1)\|y-y'\|_Y+s_1\|y-y'\|_Y))=\gamma\|y-y'\|_Y.
\end{eqnarray*}
\end{itemize}
Since we have considered all possibilities for $y,y'$, we conclude that $\Phi_0$ is $\gamma$-Lipschitz. 

We now move to proving (ii).
It follows from the definition of $\phi_j$ that it is a function  supported on $B^j$ and that $\Phi(y_j)=\phi(y_j)=\sigma_j f_j$. Next, we show that $\phi_j$ is Lipschitz. Indeed, for $y,y'\in  \bar  B^j$, we have
\be 
\label{LipOnCube}
\|\phi_j(y)-\phi_j(y')\|_X=\left|\|y_j-y\|_Y-\|y_j-y'\|_Y\right|\frac{|\sigma_j|}{\rho_j} \|f_j\|_X
\leq \frac{|\sigma_j|}{\rho_j} \|y-y'\|_Y,
\ee
where we have used that $ \|f_j\|_X=1$. Clearly, if $y,y'$ belong to the complement of $B^j$, we have
$\|\phi_j(y)-\phi_j(y')\|_X=0$. The third case is when  $y\in B^j$ and $y'$ is in its complement. Let $y_0$ be such that 
$y_0=sy+(1-s)y'$ for some $s\in [0,1]$ and $\|y_0-y_j\|_Y=\rho_j$, that is,  $y_0$ is the intersection of $\partial B^j$ and the line segment connecting $y$ with $y'$.
In this case $\phi_j(y')=\phi_j(y_0)=0$, and we can use \eref{LipOnCube},
$$
\|\phi_j(y)-\phi_j(y')\|_X=\|\phi_j(y)-\phi_j(y_0)\|_X\leq \frac{|\sigma_j|}{\rho_j} \|y-y_0\|_Y=\frac{|\sigma_j|}{\rho_j} (1-s)\|y-y'\|_Y\leq \frac{|\sigma_j|}{\rho_j}\|y-y'\|_Y.
$$
Thus, $\phi_j$ is a Lipschitz function on $Y$ with Lipschitz constant $|\sigma_j|/\rho_j$. The fact that the sum $\Phi$ is a Lipschitz mapping with the advertised Lipschitz constant follows from (i)  with $\varphi_j=\phi_j$. The proof is completed.
\hfill $\Box$

We use the  Lipschitz function $\Phi$, constructed in Lemma \ref{map} to prove the following theorem.

\begin{theorem}
	\label{TB}
	For any compact subset $\cK\subset X$ of a Banach space $X$ and   any $n\geq 1$ we have that 
	\be
	\label{stat}
	 d^{2^k\rad(\cK)}_{n}(\cK)_X\leq \e_{ kn}(\cK)_X, \quad { k=1,2,\dots.}
	\ee
In particular, when $k=n$, we have
\be
\label{power}
d^{2^n\rad(\cK)}_{n}(\cK)_X\leq \e_{ n^2}(\cK)_X, \quad { n=1,2,\dots.}
\ee	
\end{theorem}
\noindent
{\bf Proof:}  We fix $k\in N$.
Let $\cK\subset X$ be a compact set in a Banach space $X$,  let $\eta>0$,  and let 
$$ 
\cX_{ kn}:=\{f_1',\ldots,f_{2^{ kn}}'\}\subset \cK
$$
be the set such that  for every $f\in \cK$ we can find $f'_{j}\in { \cX_{kn}}$ such that 
\be
\label{no} 
\|f-f'_j\|_X\leq \e_{ kn}(\cK)_X +\eta.
\ee
 Since $\cK$ is bounded,  we can assume that $\cK\subset B_X(0,r)$ for some $r>0$.
Let us divide the  unit ball $(B_n,\|\cdot\|_{\ell_\infty^n})=[-1,1]^n\subset \R^n$   into { $2^{kn}$} non-overlapping 
 open balls $B^j$, each of side length {$2^{1-k}$}. Let us denote by $y_j$ the center of $B^j$
and  define a map $\phi_j:\R^n\to X$ as
$$
\phi_j(y)=\left(1- 2^{k}\|y_j-y\|_{\ell_\infty^n}\right)_+\cdot f_j^\prime  \in X, \quad j=1,\ldots,{ 2^{kn},}
$$
and
$$
\Phi:=\sum_{j=1}^{{2^{kn}}} \phi_j.
$$
We apply  Lemma \ref{map},  (ii) with $\sigma_j=\|f'_j\|_X$, $f_j=\frac{1}{\|f'_j\|_X}f'_j$, $\rho_j={ 2^{-k} }$,  $Y=(\R^n,\|\cdot\|_{\ell_\infty^n})$ and conclude that 
$\Phi:Y\to X$ is a map with  Lipschitz constant  
$
\gamma:=2^k\max_j\|f_j'\|_X\leq 2^kr,
$
 and $\Phi(y_j)=f_j^\prime$.
Therefore, we have   $d_n^{2^kr}(\cK)_X\leq \e_{kn}(\cK)_X+\eta$, and taking $\eta\to 0$, we obtain
\be
\label{wq1}
 d_n^{2^kr}(\cK)_X\leq \e_{kn}(\cK)_X, \quad n\geq 1.
\ee
Now, for any $\e>0$ we can find   $g=g(\e)\in X$  such that 
$$
\sup_{f\in \cK}\|f-g\|_X<\rad(\cK)+ \e2^{-k},
$$
We apply \eref{wq1} for the set $(\cK-g)$  with $r=\rad(\cK)+ \e2^{-k}$ and using  Remark \ref{zeroremark}, we arrive at 
$$
d^{2^k\rad(\cK)+\e}_{n}(\cK)_X\leq \e_{kn}(\cK)_X.
$$
The statement \eref{stat} of the theorem is obtained from the latter inequality using the continuity of the Lipschitz width $d^{\gamma}_{n}(\cK)_X$ with respect to $\gamma$, see Theorem \ref{continuous}.
\hfill $\Box$

We want to point out that estimate	\eref{stat}  in Theorem \ref{TB} is almost optimal as the following example shows.

\noindent
  \begin{example} \rm
   We consider the Hilbert space $H$ which we identify with the sequence space 
  $$
   \ell_2:=\{x=(x_1,\ldots,x_j,\ldots):\,\,\|x\|_{\ell_2}^2=\|x\|_{H}^2=\sum_{j=1}^\infty x_j^2<\infty, \,\,x_j\in \R\}.
  $$
   For each $n=1,2,\dots$,  we construct  the  compact set $\cK_n$,
  $$
\cK_n:=\{e_1,e_2,\dots, e_{2^n}, e_{2^n+1}\}\subset \ell_2,
$$
where $(e_j)$ is the standard basis in $\ell_2$, that is, all coordinate components of $e_j$ are $0$'s, except the $j$-th, which is $1$.  Then we have
$$
\tilde\e_n(\cK_n)_H\leq 3 d_{n/7}^{2\diam(\cK_n)} (\cK_n)_H.
$$
Indeed, since $\|e_i-e_j\|_{H}=\sqrt{2}$, $i\neq j$, it follows that 
$\tilde\e_k(\cK_n)_H=\sqrt 2$, $k\leq n$.
Now suppose that we have $d^\gamma_s(\cK_n)_H< \sqrt 2/3$ for some $s$ and 
$\gamma$. This means that there exists 
a norm $\|\cdot\|_{Y_s}$ on $\R^s$ and a $\gamma$-Lipschitz map $\phi$,
$\phi:(B_{Y_s},\|\cdot\|_{Y_s})\rightarrow  H$, defined on the unit ball $B_{Y_s}$, with the property that $\|\phi(y^j)-e_j\|_H< \sqrt 2/3$, $j=1,\ldots,2^n+1$, $y^j\in B_{Y_s}$, 
$j=1,\ldots,2^n+1$.
Since for $i\neq j$,
$$
\gamma\|y^j-y^{i}\|_{Y_s}\geq \|\phi(y^j)-\phi(y^{i})\|_H=\|(\phi(y^j)-e_j )+(e_j-e_{i})+(e_{i} -\phi(y^{i}))\|_H> \sqrt 2-2\sqrt 2/3 =\sqrt 2/3,
$$
we have that  $\{y^j\}_{j=1}^{2^n+1}$ is $\sqrt 2/(3\gamma)$ packing of $B_{Y_s}$. 
Using \eref{known}, we obtain that 
$$
\widetilde \cN_{\sqrt 2/(6\gamma)}(B_{Y_s})\geq 2^n+1\quad \Rightarrow \quad  \tilde \e_n(B_{Y_s})_H\geq \sqrt 2/(6\gamma).
$$
On the other hand, it follows from   \cite{CS} that 
 $4\cdot 2^{-n/s}\geq  \e_n(B_{Y_s})_H\geq 2^{-1}\tilde \e_n(B_{Y_s})_H$, and therefore,
$$
\sqrt2/(12\gamma)\leq 4 \cdot 2^{-n/s}.
$$
 Thus, for any pair $(\gamma,s)$ such that $\sqrt2/(12\gamma)> 4 \cdot 2^{-n/s}$ we get   
 $ 3 d_s^\gamma (\cK_n)_H\geq \sqrt 2=\tilde\e_n(\cK_n)_H$. This holds, for example, when $\gamma=2\diam (\cK_n)=2\sqrt 2$ and  $s= n/7$.
\end{example}

\begin{cor} 
\label{tend0}
 For every  compact subset  $\cK\subset X$ of a Banach space $X$ and every   $\gamma \geq 2\rad(\cK)$  we have
$$
\lim_{n\to \infty} d_n^{\gamma}(\cK)_X=0
$$ 
\end{cor}
\noindent
{\bf Proof:} This follows from Theorem \ref{TB},  Lemma \ref{LL1}, (iv)  and the fact that  $\lim_{n\to\infty}\e_n(\cK)_X=0$ for compact sets $\cK$,  see Remark \ref{totally bounded}.

\subsection{Estimates for Lipschitz widths from below}

We start this section with 
 a lower bound on the  Lipschitz constant $\gamma$ in $d_n^\gamma(\cK)_X$. The following proposition holds.

\begin{prop}
\label{carl}
If $d^\gamma_n(\cK)_X<\e$ for a compact subset $\cK\subset X$ of a Banach space $X$, then 
\be
\label{sas1}
\gamma\geq \frac{1}{3}\e N_{2\e}^{1/n}(\cK),
\ee
 where $N_\e(\cK)$ is the $\e$-covering number of $\cK$. In particular, if $d^\gamma_n(B_{Z_m})_X<\e$, then 
\be
\label{sasa}
 \gamma\geq \frac{1}{3}2^{-m/n}\e^{1-m/n}.
\ee
\end{prop}
\noindent
{\bf Proof:} 
 If  $d^\gamma_n(\cK)<\e$, then there is a $\gamma$-Lipschitz map $\Phi$ and a norm $\|\cdot\|_{Y_n}$, $\Phi:(B_{Y_n},\|\cdot\|_{Y_n})\to { X}$ such that $\Phi(B_{Y_n})$ approximates $\cK$ up to accuracy $\e$.	
Let us consider $\Phi(B_{Y_n})$ and let $\{y_j\}_{j=1}^N\subset  B_{Y_n}$  be such that $\{\Phi(y_j)\}_{j=1}^N$ is a maximal $\e$-packing of $\Phi(B_{Y_n})$. Then, we have 
		$$
		\e<\|\Phi(y_j)-\Phi(y_{j'})\|_X\leq \gamma\|y_j-y_{j'}\|_{Y_n},
		$$
		and thus 
		$$
		\|y_j-y_{j'}\|_{Y_n}>\e\gamma^{-1}, \quad j\neq j', \quad j,j'=1,\ldots,N.
		$$
		Therefore, see e.g. \cite[Chp. 15 Prop. 1.3]{LGM},
		\be
		\label{lal}
		N\leq  \widetilde P_{\e\gamma^{-1}}(B_{Y_n})\leq 3^n(\e\gamma^{-1})^{-n}=\left(\frac{3}{\e}\right)^n\gamma^n.
		\ee
		For every $z\in \cK$ we can find 
		$\Phi(y)$, $y\in B_{Y_n}$ such that $\|z-\Phi(y)\|_X< \e$ since $\Phi(B_{Y_n})$ approximates $\cK$ up to accuracy $\e$.  
		Since the set $\{\Phi(y),\Phi(y_1),\ldots,\Phi(y_N)\}$ is not an $\e$-packing  for $\Phi(B_{Y_n})$,  there is index $j_0$, $1\leq j_0\leq N$, 
		such that $\|\Phi(y)-\Phi(y_{j_0})\|_X\leq \e$. Then, 
		$$
		\|z-\Phi(y_{j_0})\|_X\leq \|z-\Phi(y)\|_X+\|\Phi(y)-\Phi(y_{j_0})\|_X < 2\e,
		$$
		and thus $\{\Phi(y_j)\}_{j=1}^N$	is a $2\e$-covering of $\cK$, which gives 
	         $$
		N\geq N_{2\e}(\cK).
		$$
		Combining the latter estimate with \eref{lal} gives \eref{sas1}.
		In particular, when $\cK=B_{Z_m}$, we know that 		
		$$
		 N_{2\e}(B_{Z_m})\geq (2\e)^{-m},
		$$
		and therefore we obtain \eref{sasa}.
The proof is completed.
		\hfill $\Box$

\begin{lemma}
\label{L3}
Let $\cK\subset X$ be a compact set and $\gamma>0$ be a fixed constant. If there is 
 $n>n_0$, $n_0=n_0(c_0,\alpha,\beta)$ such that 
$$d^\gamma_n(\cK)_X< c_0\frac{[\log_2 n]^{\beta}}{n^{\alpha}},  \quad \hbox{with}\quad  \alpha>0,\quad \hbox{and} \quad \beta\in \R,
$$
 then
\be
\label{sa}
 \e_m(\cK)_X < C\frac{[\log_2 m]^{\alpha+\beta}}{m^{\alpha} },\quad \hbox{with}\,\,m=cn\log_2n,
\ee
where $C,c$ are fixed constants, depending only on $\gamma$, $c_0$,  $\alpha$ and $\beta$.
\end{lemma}
\noindent
{\bf Proof:} 
 We use Proposition \ref{carl} with $\e=c_0[\log_2 n]^{\beta}n^{-\alpha}$ to obtain that 
$$
N_{2\e}(\cK)\leq \left(\frac{3\gamma}{\e}\right)^n=(3\gamma c_0^{-1}[\log_2 n]^{-\beta}n^{\alpha})^n < 
2^{n(\log_2 (3\gamma c_0^{-1})+\alpha\log_2 n-\beta\log_2(\log_2 n))}
<2^{cn\log_2 n},
$$
and therefore
$$
\e_{cn\log_2n}(\cK)_X   \leq 2c_0 [\log_2 n]^{\beta}n^{-\alpha}.
$$
If we set $m=cn\log_2n>cn$, then 
 $n=m/c\log_2 n$ and we get
\be
\label{d1}
\e_m(\cK)_X\leq  2c_0[\log_2 n]^\beta [m/c\log_2 n]^{-\alpha}=
 2c_0c^\alpha  m^{-\alpha} [\log_2 n]^{\beta+\alpha}.
\ee
Since
$$
\log_2 m=\log_2 c +\log_2 n +\log_2 \log_2 n,
$$
 for $n$ sufficiently  big we have
$$
2^{-1}\log_2 n<\log_2m < 3\log_2 n,
$$ 
and the statement follows from \eref{d1}.
\hfill $\Box$

Lemma {\rm \ref{L3}} is similar to the classical Carl's inequalities \cite{C1}, traditionally used to provide lower bounds. However,  there is an important difference. Note that Lemma 
{\rm \ref{L3}}  works for each $n$ separately,  whenever the Carl's inequality requires an assumption for   all  $j\leq n$. On the other hand the Carl's  inequality gives the upper bound for $\e_n$ not $\e_m$.

Next, we continue with a series of results presenting lower   bounds for the Lipschitz widths of compact sets, provided we have information about the entropy numbers of these sets. 
We start with a natural consequence of Proposition \ref{carl}.

\begin{prop}\label{carl2} 
Let  $\cK\subset X$ be a compact set and let 
$$
\e_n(\cK)_X> \eta_n, \quad  n=1,2,\dots,
$$
 where  $(\eta_n)_{n=1}^\infty$ is a sequences of real numbers decreasing to zero.  
Let for some $m\in \N$ { and some $\delta>0$}
$$
 d_{m}^\gamma(\cK)_X<\delta.
$$
 Then we have
\be
\label{carl3}
{ \eta_{m\log_2 (3\gamma\delta^{-1})} { <}2 \delta.}
\ee
\end{prop}
\noindent
 {\bf Proof:}
We apply Proposition \ref{carl} with $\e=\delta$ and obtain
$$
N_{2\delta}(\cK)\leq \left(\frac{3\gamma}{\delta}\right)^{ m}=2^{{m}\log_2 (3\gamma\delta^{-1})}.
$$ 
Using our assumptions and the definition of entropy numbers,  we derive
$$
2\delta\geq \e_{{ m}\log_2 (3\gamma\delta^{-1})}(\cK)_X> \eta_{{ m}\log_2 (3\gamma\delta^{-1})}. 
$$
\hfill $\Box$

 The next theorem discusses lower bounds of the Lipschitz widths $d_n^\gamma(\cK)_X$ in the case when $\gamma>0$ is a fixed constant.

\begin{theorem}
\label{widthsfrombelow}  
For any  compact set $\cK\subset X$ the following holds:
\begin{enumerate}
\item  If for some constants $c_1>0,\alpha>0$ and $\beta\in \R$ we have 
$$
\e_n(\cK)_X> c_1 \frac{(\log_2 n)^\beta}{n^{\alpha}}, \quad n=1,2,\dots,
$$
 then for each $\gamma>0$ there exists a constant  $C>0$ such that
\be 
\label{widths(i)}
d_n^\gamma(\cK)_X\geq C\frac{(\log_2 n)^{\beta-\alpha}}{n^\alpha}, \quad n=1,2,\dots.
\ee

\item  If for some constants $c_1>0,\alpha>0$  we have  
$$
\e_n(\cK)_X> c_1 (\log_2 n)^{-\alpha}, \quad n=1,2,\dots,
$$ 
then for each $\gamma>0$ there exists a constant  $C$ such that 
\be 
\label{widths(ii)}
d_n^\gamma(\cK)_X\geq C(\log_2 n)^{-\alpha}, \quad n=1,2,\dots .
\ee
 
\item 
If for some constants $c_1, c>0$  and $1>\alpha>0$  we have  
$$
\e_n(\cK)_X  > c_1 2^{-cn^\alpha}, \quad n=1,2,\dots,
$$
 then for each  $\gamma\geq 2 \rad(\cK)$ we have 
\be 
d_n^\gamma(\cK)_X\geq C 2^{-c_2  n^{\alpha/(1-\alpha)} } , \quad n=1,2,\dots,
\ee
where $C,c_2 >0 $ are constants depending on $\gamma$, $c$, and $\alpha.$
\end{enumerate}
\end{theorem}
\noindent
{\bf Proof:} 
We prove  (i)  by  contradiction. 
 If \eref{widths(i)} does not hold for some constant $C$,  then
there exists a strictly  increasing sequence of integers $(n_k)_{k=1}^\infty$, such that
$$ 
a_k:= \frac{d_{n_k}^\gamma(\cK)_Xn_k^\alpha}{(\log_2 n_k)^{\beta-\alpha}}\to 0\quad\hbox{as}\quad k\to \infty.
$$
Thus, we can write
\be 
d_{n_k}^\gamma (\cK)_X= \frac{a_k\left[\log_2 n_k\right]^{\beta-\alpha}}{n_k^\alpha}<
\frac{2a_k\left[\log_2 n_k\right]^{\beta-\alpha}}{n_k^\alpha}=:\delta_k \mbox{ for } k=1,2,\dots.
\ee
Now  we apply Proposition \ref{carl2} with $\eta_n=c_1 \frac{(\log_2 n)^\beta}{n^{\alpha}}$
and obtain
$$
c_1 \left[\log_2 ({{n_k}} \log_2(3\gamma\delta_{k}^{-1}))\right ]^\beta  {{{n_k}}}^{-\alpha} \left[\log_2(3\gamma\delta_{k}^{-1})\right]^{-\alpha } \leq 4 \frac{
{ a_k}\left[\log_2 {n_k}\right]^{\beta-\alpha}}{ {n_k}^\alpha},
$$
which we rewrite as
\be
\label{f1}
\left[\log_2  {n_k} + \log_2( \log_2(3\gamma\delta_{k}^{-1}))\right]^\beta  
\left[\log_2(3\gamma\delta_{k}^{-1})\right]^{-\alpha } \leq C_1 a_k \left[\log_2 {n_k}\right]^{\beta-\alpha},\quad \mbox{ where } C_1=4/c_1.
\ee
Observe that
$$  
\log_2(3\gamma\delta_{k}^{-1})=
\log_2 {(1.5\gamma)} +\log_2 a_k^{-1} +\alpha \log_2 {n_k} +(\alpha-\beta) \log_2 (\log_2 {n_k}),
$$ 
and therefore  for $k$ big enough we obtain
\be
\label{gq}
 \log_2 (3\gamma\delta_{k}^{-1})\leq 2\left[\log_2 (a_k^{-1}) +\alpha \log_2 n_k\right ].
\ee
The latter inequality and  \eref{f1} give
$$
2^{-\alpha}\left[\log_2  {n_k} + \log_2 (\log_2(3\gamma\delta_{k}^{-1}))\right]^\beta \left[\log_2 (a_k^{-1}) +\alpha \log_2 n_k\right ]^{-\alpha}
\leq C_1 a_k \left[\log_2 {n_k}\right]^{\beta-\alpha},
$$
which is equivalent to
\be
\label{f2}
a_k^{-1}\left[\log_2  {n_k} + \log_2 (\log_2(3\gamma\delta_{k}^{-1}))\right]^\beta
\leq 2^\alpha C_1\left[\frac{\log_2 (a_k^{-1})}{\log_2 n_k }+\alpha \right ]^{\alpha}\left[\log_2 {n_k}\right]^{\beta}.
\ee
Note that since $\delta_k\to 0$ as $k\to 0$, we have that for $k$ big enough $\log_2( \log_2(3\gamma\delta_{k}^{-1}))>0$. Now we consider several cases.

\noindent
{\bf Case 1:}   $\beta\geq 0$. In this case we have for $k$ big enough
$$
\left[\log_2  {n_k}\right]^\beta\leq \left[\log_2  {n_k} + \log_2( \log_2(3\gamma\delta_{k}^{-1}))\right]^\beta  
$$
and therefore it follows from \eref{f2} that
$$
a_k^{-1}\leq  2^\alpha C_1\left[\frac{\log_2 (a_k^{-1})}{\log_2 n_k }+\alpha \right ]^{\alpha}<C[\log_2 (a_k^{-1})]^{\alpha},
$$
which contradicts the fact that $a_{k}\to 0$ (and thus $a_k^{-1}\to \infty$).

\noindent
{\bf Case 2:}   $\beta< 0$. In this case we have 
\be
\label{kq}
\left[\log_2  {n_k}+\log_2(3\gamma\delta_{k}^{-1})\right]^\beta<\left[\log_2  {n_k} + \log_2 (\log_2(3\gamma\delta_{k}^{-1}))\right]^\beta, 
\ee
and therefore it follows from \eref{f2} that
$$
a_k^{-1}\left[\log_2  {n_k} +  \log_2(3\gamma\delta_{k}^{-1})\right]^\beta
\leq 2^\alpha C_1\left[\frac{\log_2 (a_k^{-1})}{\log_2 n_k }+\alpha \right ]^{\alpha}\left[\log_2 {n_k}\right]^{\beta}.
 $$
This gives, using \eref{gq}
\begin{eqnarray*}
a_k^{-1}
&\leq &2^\alpha C_1\left[\frac{\log_2 (a_k^{-1})}{\log_2 n_k }+\alpha \right ]^{\alpha}
\left[1+\frac{\log_2(3\gamma\delta_{k}^{-1})}{\log_2  {n_k} }\right]^{-\beta}\\
&\leq&
2^\alpha C_1\left[\frac{\log_2 (a_k^{-1})}{\log_2 n_k }+\alpha \right ]^{\alpha}
\left[1+2\alpha+2\frac{\log_2 (a_k^{-1})}{\log_2  {n_k} }\right]^{-\beta}<C[\log_2 (a_k^{-1})]^{\alpha-\beta},
\end{eqnarray*}
 which also contradicts the fact that $a_{k}\to 0$ (and thus $a_k^{-1}\to \infty$).}

To prove (ii), we repeat the argument for (i), namely, we assume that (ii) does not hold. Therefore there exists a strictly  increasing sequence of integers 
$(n_k)_{k=1}^\infty$, such that
$$ 
b_k:= d_{n_k}^\gamma(\cK)_X[\log_2 n_k]^\alpha\to 0\quad\hbox{as}\quad k\to \infty.
$$
We write
\be 
d_{n_k}^\gamma (\cK)_X= b_k[\log_2 n_k]^{-\alpha}<
2b_k[\log_2 n_k]^{-\alpha}=:\delta_k \mbox{ for } k=1,2,\dots,
\ee
and use  Proposition \ref{carl2}  with $\eta_n=c_1 (\log_2 n)^{-\alpha}$
to derive
$$
c_1 \left[\log_2 (n_k \log_2(3\gamma\delta_{k}^{-1}))\right ]^{-\alpha}
 \leq 4b_k[\log_2 {n_k}]^{-\alpha}.
$$
The latter inequality is equivalent to
\be
\label{f3}
[\log_2  {n_k} + \log_2 (\log_2(3\gamma\delta_{k}^{-1}))]^{-\alpha } \leq C_1 b_k (\log_2 {n_k})^{-\alpha},\,\, C_1:=4/c_1,
\ee
which, after using \eref{kq}  with $\beta=-\alpha$ gives
$$
[\log_2  {n_k} +  \log_2(3\gamma\delta_{k}^{-1})]^{-\alpha } \leq C_1 b_k \left[\log_2 {n_k}\right]^{-\alpha}.
$$
We continue by writing the above inequality as
$$
b_k^{-1}\leq C_1 \left[1+  \frac{\log_2(3\gamma\delta_{k}^{-1})}{\log_2  {n_k} }\right]^{\alpha } \leq
\left[1+2\alpha+2\frac{\log_2 (b_k^{-1})}{\log_2  {n_k} }\right]^{\alpha }\leq C[\log_2 (b_k^{-1})]^{\alpha},
$$
where we have used \eref{gq}.
The latter inequality  contradicts the fact that $b_{k}$ tends to zero, and the proof of (ii) is completed.

We now  prove (iii). To simplify the notation, we denote by $d_n:=d_n^{\gamma} (\cK)_X$   and observe that,
according to Corollary \ref{tend0}, $d_n\to 0$ for $n\to\infty$ when $\gamma\geq \rad(\cK)$. We use Proposition \ref{carl2} with $\delta=2d_n$ and  $\eta_n=c_12^{-cn^\alpha}$ to obtain the inequality
\be \label{exp1}
4d_n \geq c_1 2^{-c[n\log_2 (3\gamma d_n^{-1})]^\alpha},
\ee
which can be rewritten  as
$$
2^{c[n\log_2 (3\gamma d_n^{-1})]^\alpha}\geq \frac{c_1}{4}d_n^{-1}\quad\Leftrightarrow\quad 2^{c[n(\log_2 \xi_n)]^\alpha}\geq 
\frac{c_1}{12\gamma}\xi_n=:A\xi_n, \,\, \hbox{where}\,\,\xi_n:=3\gamma d_n^{-1}\to \infty \,\,\hbox{as}\,\,n\to\infty.
$$
Taking  logarithm on both sides of the inequality  and using the fact that $\xi_n\to  \infty$ we obtain for $n$ big enough
\be\label{exp3}
\log_2 \xi_n\leq c n^\alpha (\log_2 \xi_n)^\alpha -\log_2 A\leq 2cn^\alpha (\log_2 \xi_n)^\alpha,
\ee 
and therefore $ \log_2\xi_n\leq (2c)^{1/(1-\alpha)} n^{\alpha/(1-\alpha)}$. Returning back to the notation for the Lipschitz width, we obtain
$$
3\gamma 2^{-(2c)^{1/(1-\alpha)} n^{\alpha/(1-\alpha)} }\leq d_n^\gamma(\cK)_X
$$
for $n$ big enough. 
This completes the proof of (iii) by choosing the constants appropriately so that the above inequality holds for all $n$.
\hfill $\Box$

\subsection{Summary}
Now we are ready to state a corollary to Theorem \ref{TB} and Theorem \ref{widthsfrombelow}.

\begin{cor}
\label{cr1}
Let  $\cK\subset X$ be a  compact subset of a Banach space $X$,  $n\in \N$, and $d_n^{\gamma}(\cK)_X$ be the   Lipschitz width for $\cK$ with Lipschitz constant  
$\gamma\geq 2\rad(\cK)$. Then the following holds:
\begin{enumerate}
\item For  $\alpha>0$, $\beta\in \R$, we have 
\begin{eqnarray}
\nonumber
\e_n(\cK)_X\leq C \frac{[\log_2n]^\beta}{n^\alpha}, \quad n=1,2,\ldots,\quad  \Rightarrow\quad 
d_n^{\gamma}(\cK)_X\leq C \frac{[\log_2n]^\beta}{n^{\alpha}}, \quad n=1,2,\ldots,\\
\nonumber
\e_n(\cK)_X\geq C \frac{[\log_2n]^\beta}{n^\alpha},\quad n=1,2,\ldots,\quad \quad \Rightarrow\quad 
d_n^{\gamma}(\cK)_X\geq C'\frac{[\log_2n]^\beta}{n^{\alpha}[\log_2 n]^{\alpha}},\quad, n=1,2,\ldots .
\end{eqnarray}
\item For  $\alpha>0$,  we have 
\be
\label{c2}
\e_n(\cK)_X\asymp \frac{1}{[\log_2 n]^\alpha},\quad n=1,2,\ldots,\quad \quad \Rightarrow\quad 
d_n^{\gamma}(\cK)_X\asymp \frac{1}{[\log_2 n]^{\alpha}},\quad n=1,2,\ldots\quad.
\ee
\item For  $0<\alpha<1$,  we have 
\begin{eqnarray}
\nonumber
\e_n(\cK)_X\leq C 2^{-cn^{\alpha}}, \quad n=1,2,\ldots,\quad  \Rightarrow\quad 
d_n^{\gamma}(\cK)_X\leq C 2^{-cn^{\alpha}}, \quad n=1,2,\ldots,\\
\nonumber
\e_n(\cK)_X\geq C 2^{-cn^{\alpha}},\quad n=1,2,\ldots,\quad \quad \Rightarrow\quad 
d_n^{\gamma}(\cK)_X\geq C'2^{-c'n^{\alpha/(1-\alpha)}},\quad n=1,2,\ldots.
\end{eqnarray}

\end{enumerate}
\end{cor}
\noindent
{\bf Proof:} We first prove (i). We  assume that 
$$
\e_n(\cK)_X\leq  C\frac{[\log_2n]^\beta}{n^\alpha},  
$$
and use  \eref{stat} in Theorem \ref{TB} with $k=1$ to derive
$$
d_n^{2\rad(\cK)}(\cK)_X\leq C\frac{[\log_2n]^\beta}{n^{\alpha}}.
$$
It follows from 
Lemma \ref{LL1}, (iv) that for $\gamma\geq 2\rad(\cK)$, we have
$
d_n^{\gamma}(\cK)_X\leq d_n^{2\rad(\cK)}(\cK)_X,
$
and combining the above two inequalities gives for all $n$
$$
d_n^{\gamma}(\cK)_X\leq C\frac{[\log_2n]^\beta}{n^{\alpha}}.
$$
The other direction in \eref{c1} is the statement of Theorem \ref{widthsfrombelow}, (i).
The proofs of (ii) and (iii) are similar and we omit them.
\hfill $\Box$

\subsubsection{Lipschitz widths   could  be smaller than entropy}
{ In this section, we show that  the estimates  in Corollary \ref{cr1} are sharp and cannot be improved.

\subsubsection{The logarithm in Corollary \ref{cr1}, part (i) cannot be removed}
Here, we provide an example of a compact set $\cK$ for which the the entropy numbers behave like 
$ n^{-1}$, while the Lipschitz width behaves as  $[n\log_2 (n+1)]^{-1}$.

We consider the Banach space $X= {\bf c_0}$ of all sequences that converge to $0$, equipped with the $\ell_\infty$ norm and its compact subset 
\be
\label{ds}
\cK(\sigma):=\{\sigma_j e_j\}_{j=1}^\infty \cup \{0\}\subset {\bf c_0},
\ee
determined by the strictly decreasing converging to $0$ sequence $\sigma:=(\sigma_j)_{j=1}^\infty$, where $(e_j)_{j=1}^\infty$ are  the standard  basis in  $c_0$.
Since
$$
\|\sigma_j e_j -\sigma_{j^\prime}e_{j^\prime}\|_{\ell_\infty}= \sigma_j, \quad \hbox{for all}\quad j'>j,
$$
it follows that the ball with center $\sigma_j e_j$ and radius $\sigma_j$  contains all points $\sigma_{j^\prime} e_{j^\prime}$ with $j^\prime>j$ and none with $j^\prime <j$. Thus, if we look for $2^n$ balls with centers in $\cK(\sigma)$ covering $\cK(\sigma)$ with smallest radius, we take the balls $B(\sigma_j e_j, \sigma_{2^n})$,  $j=1,2,\dots, 2^n$, with centers $\sigma_je_j$ and radius $\sigma_{2^n}$. Each of the first $2^n-1$ balls contain only one point from $\cK(\sigma)$, while the last ball $B(\sigma_{2^n} e_{2^n}, \sigma_{2^n})$ contains the rest of the points $\{\sigma_{j} e_j\}_{j=2^n}^\infty \cup \{0\}$, which gives 
\be
\label{entrK}
\tilde \e_n(\cK(\sigma))_{X}=\sigma_{2^n}.
\ee

 We next investigate the behavior of  $d_n^\gamma(\cK(\sigma))_X$. We shall use the following lemma  which gives upper bounds for the Lipschitz widths for the sets
 $\cK(\sigma)$.

\begin{lemma}
\label{refr}
 Consider the strictly decreasing sequence $\sigma:=(\sigma_j)_{j=1}^\infty$, $\sigma_j\to 0$ as $j\to\infty$,  and the set 
$\cK(\sigma)$, defined in \eref{ds}. 
If  $\sigma_1\leq \gamma/2$ and we can find  $N$ (finite or infinite) such that 
\be 
\label{volumfit}
\sum_{j=1}^N \sigma_j^n\leq (\gamma/2)^n,
\ee
then $d_n^\gamma(\cK(\sigma))_X\leq \sigma_N$.
\end{lemma}
\noindent
{\bf Proof:} We consider the case when $N$ is finite. Similar arguments hold in the infinite case. 
 For every $\sigma_j$, $j=1,\ldots,N$, we define $\ell_j\in \N\cup\{0\}$ as
\be
\label{pl}
2^{-\ell_j-1}<2\frac{\sigma_j}{\gamma}\leq 2^{-\ell_j}.
\ee
Then it follows from \eref{volumfit} that 
\be
\label{bb}
\sum_{j=1}^N2^{-n\ell_j}\leq \sum_{j=1}^N \left( 4\sigma_j/\gamma\right)^n  \leq 2^n.
\ee
Since $(\sigma_j)_{j=1}^\infty$ is a decreasing sequence, we have that 
$2^{-\ell_1}\geq 2^{-\ell_2}\geq 2^{-\ell_3}\geq \dots\geq 2^{-\ell_N}$. Note that  some of the $\ell_j$'s can be equal to each other. Let $k_1,k_2,\ldots, k_s=N$, be the indices such that 
$$
\ell_1=\ldots =\ell_{k_1}<\ell_{k_1+1}=\ldots=\ell_{k_2}<\ell_{k_2+1}=\ldots=\ell_{k_s}=\ell_N.
$$
We set $k_0=0$ and   rewrite   inequality  \eref{bb} as
\be
\label{ss1}
2^n\geq \sum_{j=1}^N 2^{-n\ell_j}=\sum_{j=1}^s (k_j-k_{j-1})2^{-n\ell_{k_j}}.
\ee
Observe that  the volume of a cube with side length $2^{-\ell_{k_j}}$ is $2^{-n\ell_{k_j}}$, while the volume of $[-1,1]^n$ is $2^n$. It follows from simple volumetric considerations, that 
we  can divide naturally the cube $[-1,1]^n$ into $k_1$ open non-overlaping cubes each with side length $2^{-\ell_{k_1}}$, $(k_2-k_1)$  open non-overlaping cubes each with side 
length $2^{-\ell_{k_2}},\ldots$, $(k_s-k_{s-1})$  open non-overlaping cubes each with side length $2^{-\ell_{k_s}}$, since, according to \eref{ss1}, the sum of the  total volumes of these cubes 
does not exceed the total volume of $[-1,1]^n$.
Thus,
there exists a sequence of non-overlapping open cubes $B^j$, 
$$
B^j :=B^j(y_j,2^{-\ell_j-1})\subset  (B_{\ell_\infty^n},\|\cdot\|_{\ell_\infty^n}):=[-1,1]^n, \quad j=1,\ldots,N,
$$
 with side length $2^{-\ell_j}$.
Then, according to Lemma \ref{map}, the mapping $\Phi:({ B_{\ell_\infty^n}},\|\cdot\|_{\ell_\infty^n})\to {\bf c_0}$, defined as
$$
\Phi(y):=\sum_{j=1}^N\sigma_j(1-2^{\ell_j+1}\|y_j-y\|_{\ell_\infty^n})_+\cdot e_j
$$
is a Lipschitz mapping. Its Lipschitz constant  is $\displaystyle{\sup_{j=1,\ldots,N} \{2^{\ell_j+1}\sigma_j}\}$ and  
$ \Phi(y_j)=\sigma_je_j$, $j=1,\ldots,N$.  It follows from \eref{pl} that 
$$
\sup_{j=1,\ldots,N} 2^{\ell_j+1}\sigma_j\leq \gamma,
$$
and therefore $\Phi$ is a $\gamma$-Lipschitz mapping. On the other hand, since
$$
\sup_{j'\geq 1}\inf_{y\in B_n}\|\sigma_{j'}e_{j'}- \Phi(y)\|_{\ell_\infty}\leq \sup_{j'\geq 1}\inf_{j=1,\ldots,N}\|\sigma_{j'}e_{j'}- \Phi(y_j)\|_{\ell_\infty}=\sup_{j'\geq 1}\inf_{j=1,\ldots,N}
\|\sigma_{j'}e_{j'}-\sigma_je_j\|_{\ell_\infty}=\sigma_N,
$$
 and
$$
\inf_{y\in B_n}\|0- \Phi(y)\|_{\ell_\infty}\leq \inf_{j=1,\ldots,N}\|\Phi(y_j)\|_{\ell_\infty}=\inf_{j=1,\ldots,N}
\|\sigma_je_j\|_{\ell_\infty}=\sigma_N,
$$
it follows that $d_n^\gamma(\cK(\sigma))_X\leq \sigma_N$ and 
the proof is completed.
\hfill $\Box$

Now, we are ready to state the main theorem in this section.
\begin{theorem}
\label{mainP}
		The compact set $\cK(\sigma)\subset  {\bf c_0}$, defined in \eref{ds}, with   $\sigma=(\sigma_j)_{j=1}^\infty$ being  the sequence 
		$\sigma_j=1/\log_2(j+1)$  has  inner entropy numbers
		$$
		\tilde\e_n(\cK(\sigma))\asymp \frac{1}{n}, 
		$$
		 and Lipschitz width  
		$$
		 { d_n^\gamma (\cK(\sigma))\asymp \frac{1}{n\log_2 (n+1)}}
		$$
for any for $\gamma> 2$.		
\end{theorem}
\noindent
	{\bf Proof:} The behavior of the entropy  follows from \eref{entrK} and  the estimate from below for the Lipschitz widths follows from Theorem, \ref{widthsfrombelow}, (i). We are only left to prove the upper estimate for the width. If we 
	 show   that \eref{volumfit}
	  holds for the choice of $\sigma_j=[\log_2(j+1)]^{-1}$, $j=1,2,\ldots$,  and 
	$N=(n+1)^n$, where  $\gamma>2$ and $n$ is sufficiently large, 
	since  $\sigma_1=1\leq \gamma/2$, we  can  use Lemma \ref{refr}  to conclude that 
$d^\gamma_n(\cK(\sigma))\leq  \sigma_N=(n\log_2 (n+1))^{-1}$, for $n\geq n_0$, depending only on $\gamma$. This could conclude the proof.

 We now concentrate on proving \eref{volumfit} with $N=(n+1)^n$ for $n$ sufficiently large. We start with
	defining  $J=J(n)$ as
	$$
	 2^{J-1}\leq (n+1)^n<2^J,
	$$
	and estimate
\be 
\label{lg11}
	\sum_{j=1}^{(n+1)^n}\sigma_j^n\leq 
	\sum_{k=0}^{J-1} \sum_{j=2^k}^{2^{k+1}-1} \sigma_j^n\leq 1+\sum_{k=1}^{J-1} 2^k k^{-n}=:1+\sum_{k=1}^{J-1} q(k),	\quad \hbox{where} \,\,q(t):=2^t t^{-n}, \quad t\geq 1.
	\ee	
Simple calculation shows that  $q(t)$ is decreasing on  $[1, n/\ln 2]$ and increasing on  $[n/\ln 2,\infty)$.	
Moreover, we have that 
	\be \label{ineq_q}
	2e^{-n/t}<\frac{q(t+1)}{q(t)}=\frac{2}{\left(1+\frac{1}{t}\right)^n}< 2e^{-n/(t+1)}\leq 1/2\quad \hbox{for}\quad t\leq \frac{n}{\ln 4}-1.
	\ee
It follows from \eref{lg11} that for  $n\geq 3$,
\be
\label{s1}
\sum_{j=1}^{(n+1)^n}\sigma_j^n\leq 1+
\sum_{1\leq k{ \leq n/\ln 4}}  q(k)+\sum_{n/\ln 4  { < k} \leq n/\ln 2}
	 q(k)+\sum_{n/\ln 2< k\leq  J-1}  q(k)=:S_1(n)+S_2(n)+S_3(n).
\ee
We will provide upper bounds for each of $S_1$, $S_2$ and $S_3$.	 Clearly
$$
S_1(n)=1+\sum_{1\leq k{ \leq n}/\ln 4}  q(k)
< 1+q(1)\cdot \sum_{1\leq k{ \leq n}/\ln 4}2^{-k+1}
<1+2\cdot \sum_{k=1}^\infty 2^{ -k+1}=5,
$$
since for this range of $k$'s we have  $q(k+1)<\frac{1}{2}q(k)$,  see \eref{ineq_q}.

 Next, note that  $q$ is a decreasing  function for the range of $k$ in $S_2$,  and therefore 
	\begin{eqnarray}
	 \label{S-2}
	 \nonumber
		S_2(n)&\leq&\left( { \frac{n}{\ln 2}-\frac{n}{\ln 4}}\right)\cdot 2^{n/\ln 4}\left(\frac{n}{\ln4}\right)^{-n}<n \left(2^{1/\ln 4}\ln 4\right)^n n^{-n}<n\left(\frac{2.3}{n}\right)^n,
	\end{eqnarray}
since { $(n/\ln 2-n/\ln 4)=n/(2\ln 2) < n$ and }  $2^{1/\ln 4}\ln 4<2.3$.  For $n\geq 5$
$$
2.3n^{1/n}<2.3(1+\frac{1}{\sqrt{2}})<n \quad \Rightarrow\quad  S_2(n)<n\left(\frac{2.3}{n}\right)^n =    n\left(\frac{2.3 n^{1/n}}{n\cdot n^{1/n}}\right)^n   <1.
$$
So, we obtain
$$
S_2(n)<1\quad \hbox{for}\quad n\geq 5.
$$
To estimate $S_3$, we notice that  the biggest summand is the last one,
	\begin{eqnarray}\label{S-3}
	\nonumber
		S_3(n)&\leq &J\cdot 2^{J-1} (J-1)^{-n}<(n\log_2(n+1)+1)(n+1)^n(n\log_2(n+1)-1)^{-n}\\
		\nonumber
		&=&n\left[1+\frac{1}{n}\right]^n\left[\log_2(n+1)+\frac{1}{n}\right]\left[\log_2(n+1)-\frac{1}{n}\right]^{-n}\\
		\nonumber
		&<& 2en \left[\log_2(n+1)-\frac{1}{n}\right]\left[\log_2(n+1)-\frac{1}{n}\right]^{-n}=2en\left[\log_2(n+1)-\frac{1}{n}\right]^{1-n}.
		\nonumber
	\end{eqnarray}
Let us now consider the functions 
$$
\ell(x):=x^{1/(x-1)}, \quad r(x):=\log_2(x+1)-\frac{1}{x}.
$$
One can show that $\ell$ is a decreasing function on the interval $[5,\infty)$, while $r$ is increasing function on the same interval. Therefore, for every $n\geq 5$
$$
n^{1/(n-1)}=\ell(n)\leq \ell(5)=5^{1/4}<\log_26- \frac{1}{5}=r(5)\leq r(n)=\log_2(n+1)-\frac{1}{n},
$$
and so
$$
 n<\left(\log_2(n+1)-\frac{1}{n}\right)^{n-1}\quad \Rightarrow\quad n\left(\log_2(n+1)-\frac{1}{n}\right)^{1-n}<1.
$$
The latter inequality combined with the estimate for $S_3$ gives that 
$$
S_3(n)<2e, \quad \hbox{for}\quad n\geq 5.
$$
Finally, combining \eref{s1} with all estimates for $S_1$, $S_2$ and $S_3$, we obtain that 
$$
\sum_{j=1}^{(n+1)^n}\sigma_j^n<S_1(n)+S_2(n)+S_3(n)<6+2e\leq (\gamma/2)^n, \quad\hbox{provided}\quad n\geq \max\left\{5,\frac{\ln(6+2e)}{\ln \gamma-\ln 2}\right\}.
$$	
The proof is completed.
\hfill $\Box$

\subsubsection{ Corollary \ref{cr1}, part (iii) cannot be improved}
In this section we show that the requirement $\alpha<1$ in Theorem \ref{widthsfrombelow}, (iii) and Corollary
\ref{cr1}, (iii) is necessary. We give an example of a compact set $\cK$ with 
$\e_n(\cK)\asymp 2^{-cn}$ and Lipschitz width zero, which shows that  an estimate from below 
for $d_n^\gamma(\cK)$ in terms of $n$ is not possible.

\begin{theorem}
\label{mainP1}
		 The compact set $\cK(\sigma) \subset {\bf c_0}$, defined in {\rm\eref{ds}} with   $\sigma=(\sigma_j)_{j=1}^\infty$ being the sequence  $\sigma_j=j^{-c}$,  ${c>0}$ has inner entropy numbers
		$$
		\tilde\e_n(\cK(\sigma))\asymp 2^{ -c  n}, 
		$$
		 and Lipschitz width
		$$
		d_n^\gamma (\cK(\sigma))=0,\quad n\geq n_1,
		$$
		 for any  $\gamma> 2$, where $n_1$ depends  only on $\gamma$ and $c$.	
\end{theorem}
\noindent
{\bf Proof:} Since $c>0$ is a fixed constant, there is $n_0$ such that $nc>1$ for   $n\geq n_0+1$
 and $nc<1$ for $n\leq n_0$. Then, for every $N\geq n_0$ we have 
\begin{eqnarray*} 
\nonumber
	\sum_{j=1}^{N}\sigma_j^n& =&
	\sum_{j=1}^{n_0}  j^{-c n} +\sum_{j=n_0+1}^{N} j^{-c n} <n_0+\int_{n_0}^{\infty} x^{-cn}\,dx
	= n_0+\frac{n_0}{cn-1}\left(\frac{1}{n_0^c}\right)^n\leq
	n_0+\frac{n_0}{cn-1}\leq\left(\frac{\gamma}{2}\right)^n
\end{eqnarray*}	
for  $\gamma>2$ and $n\geq n_1(\gamma,c)\geq n_0+1$, and $\sigma_1=1< \gamma/2$. It follows from Lemma \ref{refr} that
$$
d_n^\gamma (\cK(\sigma))\leq \sigma_N=N^{-c}, \quad \hbox{for all}\quad n\geq n_1, \quad N\geq n_0.
$$
Letting $N\to\infty$ gives that $d_n^\gamma (\cK)=0$, provided $n\geq n_1$.
Finally, it is easy to show that the inner entropy numbers
$\tilde\e_n(\cK(\sigma))\asymp 2^{-c n}$, and the proof is completed.
\hfill $\Box$

}

\section{ Comparison between Lipschitz  and  Kolmogorov  widths}
\label{S5}

If we fix the value of $n\ge 0$, the Kolmogorov $n$-width of $\cK$
 is defined as
 \be
 \label{Kwidth}
  d_0(\cK)_X=\sup_{f\in \cK} \|f\|_X, \quad d_n(\cK)_X:=\inf_{\dim(X_n)=n} \sup_{f\in \cK}\dist(f,X_n) _X, \quad  n\geq 1.
 \ee
 It  tells us the optimal performance possible  for the approximation
 of the model class $\cK$ using linear spaces of dimension $n$. However, it does not tell us how to select
 a (near) optimal space $Y$ of dimension $n$ for this purpose. 
Let us note that in the definition of Kolmogorov width,  we are not requiring that   the mapping 
 which sends $f\in \cK$ into an approximation to $f$  is a linear map.    
  There is a concept of {\it linear} width which requires the linearity of the approximation map.  
Namely, given $n\ge 0$ and a model class $\cK\subset X$, its {\it linear} width   $d_n^L(K)_X$ is defined as
\be
   \label {linwidth}
     d_0^L(\cK)_X=\sup_{f\in \cK} \|f\|_X, \quad    
 d_n^L(\cK)_X:=\inf_{L\in \cL_n} \sup_{f\in \cK} \|f-L(f)\|_X, \quad  n\geq 1,
\ee
    where the infimum is taken over the class  $\cL_n$ of all  continuous linear maps from $X$ into itself with rank at most $n$.   
    
  We prove in the next theorem the   intuitive fact  that the Lipschitz width is smaller than the Kolmogorov width.

\begin{theorem}
\label{TK}
For every compact set $\cK\subset X$ and every $n\geq 1$, we have
 \be
\label{compKol}
d^\gamma_n(\cK)_X\leq d_n(\cK)_X\leq d^L_n(\cK)_X, \quad \hbox{ for } \quad \gamma=d_n(\cK)_X+\rad(\cK).
\ee 
\end{theorem}
\noindent
{\bf Proof:}  
It is clear that $d_n(\cK)_X\leq d^L_n(\cK)_X$ for every $n\geq 0$ since we can take $X_n$ to be the $n$-dimensional linear space containing $L(X)$ when $L\in \cL_n$, so we concentrate on the first inequality.
	We   start with  $\gamma>d_n(\cK)_X+\rad(\cK)$, denote 
		$$
		\eta:=\gamma-d_n(\cK)_X-\rad(\cK)>0,
		$$
		 and  choose  $\eta_1$ to be such that $0<\eta_1<\eta$.
	Let  $X_n\subset X$ be an $n$-dimensional linear subspace in $X$ such that,
	$$
	\sup_{f\in \cK} \inf_{g\in X_n}\|f-g\|_X< d_n(\cK)_X+\eta_1.
	$$
	For  every $f\in \cK$,  we denote by  $g=g(f)$ the element in $X_n$  for which
\be
\label{nnm}
	\|f-g(f)\|_X< d_n(\cK)_X+\eta_1,
\ee
and the collection of all such elements are denoted by
	$$
	\cA=\{g(f): \,f\in \cK\}\subset X_n.
	$$
	Let us fix $g_0\in X$ such that $\sup_{f\in \cK}\|f-g_0\|_X< \rad{\cK}+\eta-\eta_1$.  Then, for every $f\in \cK$, 
	$$
	\|g(f)-g_0\|_X\leq \|g(f) -f\|_X+\|f -g_0\|_X< d_n(\cK)_X+\rad (\cK) +\eta= \gamma,
	$$ 
	and therefore 
	$$
	\rad(\cA)<  \gamma, \quad \hbox{and}\quad \cA\subset B(g_0,\gamma) :=\{g\in X_n:\|g-g_0\|_X\leq \gamma\}.
	$$
	We now define the mapping  $\Phi:(B_{X_n},\|\cdot\|_X)\to X$ from the unit ball 
	$B_{X_n} :=\|g\in X_n:\,\|g\|_X\leq 1\}$ in $X_n$ as
	$\Phi(g)=g_0+ \gamma g$. Clearly $\Phi $ is a $\gamma$-Lipschitz map. Moreover, since $\Phi(B_{X_n})= B(g_0,\gamma)$ and $\cA\subset B(g_0,\gamma)$,
	we have that 
	$$	
	\sup_{f\in \cK}\inf_{g\in B_{X_n}}\|f-\Phi(g)\|_X\leq\sup_{f\in \cK}\inf_{g\in \cA}\|f-g\| < d_n(\cK)_X+\eta_1,  
        $$
	where we have used \eref{nnm} in the last inequality.
	Thus, using  Lemma \ref{LL1},  (iii),   we obtain  
	$$
	d^\gamma_n (\cK)_X\leq  d_n(\cK)_X+\eta_1,
	$$
	and letting  $\eta_1\to 0$ gives
	$$
	d^\gamma_n (\cK)_X\leq d_n(\cK)_X, \quad \hbox{for any}\,\,\gamma>d_n(\cK)_X+\rad(\cK).
	$$
	 Now  \eref{compKol} follows  from Theorem \ref{continuous}  by  taking $\gamma\to d_n(\cK)_X+\rad(\cK)$.
	\hfill$\Box$

\begin{cor}
For every $n\geq 1$ and  every compact set $\cK\subset X$ we have
\be
\label{ll}
d^\gamma_n(\cK)_X\leq d_n(\cK)_X, \quad  \gamma=2\sup_{f\in \cK}\|f\|_X.
\ee

\end{cor}
\noindent
{\bf Proof:} The  inequality  follows from   Theorem \ref{TK},  Lemma \ref{LL1}, (iv),  and the  fact that for every $n\geq 1$
$$
d_n(\cK)_X+\rad(\cK)\leq 2\sup_{f\in \cK}\|f\|_X.
$$
\hfill $\Box$

As a result of this section, we can give the following  improvement of Corollary \ref{tend0}.  
\begin{cor}
\label{tend0_1}
 If $\cK\subset X$ is compact, then for every  $n_0\in \N\cup\{0\}$   and every
 $\gamma \geq \d_{n_0}(\cK)_X+\rad(\cK)$ we have
		$$
		\lim_{n\to \infty} d_n^\gamma (\cK)_X=0.
		$$
				\end{cor}
\noindent
{\bf Proof:} The statement  follows from Theorem \ref{TK},  Lemma \ref{LL1} (iv)  and the fact that the sequence of Kolmogorov widths 
$(d_n(\cK)_X)$ of a compact set $\cK$ is a non-increasing sequence of non-negative numbers that tends to zero, see e.g.\cite[Prop 1.2]{AP}.
\hfill $\Box$

\subsection{Examples of different behavior of the Lipschitz and  Kolmogorov widths}

 It is intuitively clear that the Lipschitz widths could be much  smaller than the Kolmogorov widths. We illustrate this observation by discussing the following two examples.

\noindent
\begin{example}  \rm
This example, borrowed from Albert Cohen, arises in some partial differential equations. We denote by $\chi_a$  the characteristic function of $[a,a+1]$, $a\in[0,1]$ and
consider the univariate linear transport equation 
\be
\label{ol}
\partial_t u_a+a\partial_xu_a=0, 
\ee
with constant velocity $a\in[0,1]$ and initial condition 
\be
\label{ol2}
u_0(x)=u_a(x,0)=\chi_0(x).
\ee
We denote by 
$$
\cH:=\{\chi_a:\,a\in[0,1]\}\equiv \{u_a(x,1):\,a\in[0,1]\}
$$
the solution manifold to \eref{ol}-\eref{ol2}  evaluated at time $t=1$. We prove the following lemma for the set $\cH$.

\begin{lemma} 
\label{Albert}
The Kolmogorov width of $\cH\subset L_1[0,2]$ is
\be
\label{kk}
 (n+1)^{-1} \leq  d_n(\cH)_{L_1[0,2]}\leq 4n^{-1},  
\ee
while   its inner entropy numbers
\be
\label{ent}
\tilde \e_n(\cH)_{L_1}=2^{-n+1}.
\ee
\end{lemma}
\noindent
{\bf Proof:}
We first observe that
$
\|\chi_a-\chi_b\|_{L_1[0,2]}=2|a-b|. 
$
If we define 
$$
t_j:=(2j+1)2^{ -n-1}, \quad j=0,1,\dots, 2^n-1,
$$
to be the centers of the intervals $[j2^{ -n},(j+1)2^{ -n}]{ \subset [0,1]}$,  we have
$$
\chi_a\in B(\chi_{t_j},2^{-n+1})\quad \Leftrightarrow\quad \|\chi_a-\chi_{t_j}\|_{L_1[0,2]}\leq 2^{-n+1}
\quad \Leftrightarrow\quad  |a-t_j|\leq 2^{-n},
$$
where $B(\chi_{t_j},2^{-n+1})$ is the closed ball in $L_1[0,2]$ with center $\chi_{t_j}$ and radius $2^{-n+1}$.
 So, those balls cover $\cH$. This calculation also shows that if we have $2^n$ balls covering $\cH$ and 
 one of them has radius $<2^{-n+1}$ then some other one must have a radius $> 2^{-n+1}$. This proves  
 \eref{ent}.

To show \eref{kk}, we first observe that the $n$-dimensional space
$$
V_n:={\rm span}\{\bar \chi_j, \quad j=0,\ldots,n-1\}, 
$$
 where $\bar \chi_j$ is the characteristic function of the interval ${[2j/n,2(j+1)/n]} $
provides an error at most $4n^{-1}$ for the elements from $\cH$. Indeed,   for each $\chi_a\in \cH$, we have
$$
\|\chi_a-\sum_{j=j_1}^{j_2}\bar \chi_j\|_{L_1[0,2]}\leq 4n^{-1}, 
$$
where $j_1=j_1(a)$ and $j_2=j_2(a)$ are defined as
$$
j_1(a)= \max \{j:\,\,2j/n\leq a,  \,\,0\leq j\leq n-1\}, \quad j_2(a)=\max \{j:\,\,2j/n\leq a+1, 0\leq j\leq n-1\},
$$
and therefore $d_n(\cH)_{L_1[0,2]}\leq  4 n^{-1}$.
To prove the lower bound in \eref{kk}, we use a  well known result, see e.g. \cite[Chap. II, Prop 1.3 ]{AP},  which states that for any unit ball $U$ in a Banach space $X$ and any finite dimensional space 
$\cV_{n+1}$ of dimension $n+1$, the Kolmogorov width 
\be
\label{pop}
d_n(U\cap \cV_{n+1})_X=1.
\ee
We apply this result for the Banach space $X=L_1[0,2]$, the  unit ball $U$ in $L_1[0,2]$, and the linear space $\cV_{n+1}\subset L_p[0,1]$, defined as
$$
\cV_{n+1}=\span\{\varphi_0,\ldots,\varphi_n\},\quad \varphi_j:=\chi_{j/(n+1)}-\chi_{(j+1)/(n+1)}, \quad j=0,\ldots,n.
$$
Another representation for the $\varphi_j$'s is
$$
\varphi_j:= \chi_{[j/(n+1),(j+1)/(n+1)]} -\chi_{[1+j/(n+1),1+(j+1)/(n+1)]},  \quad j=0,\ldots,n,
$$
and  since they have disjoint supports,
every $\varphi=\sum_{j=0}^n\alpha_j\varphi_j\in \cV_{n+1}$ has norm
\be
\label{la1}
\|\varphi\|_{L_1[0,2]}=\|\sum_{j=0}^n\alpha_j\varphi_j\|_{L_1[0,2]}=2(n+1)^{-1}\sum_{j=0}^n|\alpha_j|.
\ee
Therefore
\be
\label{q111}
\varphi=\sum_{j=0}^n\alpha_j\varphi_j\in U\cap\cV_{n+1}\quad \Leftrightarrow\quad \sum_{j=0}^n|\alpha_j|\leq \frac{1}{2}(n+1).
\ee
Let us fix an $n$ dimensional subspace $V_n$ and let $v_j\in V_n$ be such that 
$$
\dist(\chi_{j/(n+1)}, V_n)_{L_1[0,2]}=\|\chi_{j/(n+1)}-v_j\|_{L_1[0,2]}, \quad j=0,\ldots,n.
$$
Then, for every $\varphi\in U\cap\cV_{n+1}$,  we have 
\begin{eqnarray*}
\dist(\varphi,V_n)_{L_1[0,2]}&\leq &\|\sum_{j=0}^n\alpha_j\varphi_j-\sum_{j=0}^n\alpha_j(v_j-v_{j+1})\|_{L_1[0,2]}\\
&=&
\|\sum_{j=0}^n\alpha_j(\chi_{j/(n+1)}-v_j)-\sum_{j=0}^n\alpha_j(\chi_{(j+1)/(n+1)}-v_{j+1})\|_{L_1[0,2]}\\
&\leq&2\dist(\cH,V)_{L_1[0,2]}\sum_{j=0}^n|\alpha_j|\leq (n+1)\dist(\cH,V_n)_{L_1[0,2]},
\end{eqnarray*}
where we have used  \eref{q111}.
Therefore, it follows from  \eref{pop} and the latter estimate that
\begin{eqnarray*}
1&=&d_n(U\cap \cV_{n+1})_{L_1[0,2]}=\inf_{V_n}\sup_{\varphi\in U\cap\cV_{n+1}}\dist(\varphi,V_n)_{L_1[0,2]}\\
&\leq& (n+1)\inf_{V_n}\dist(\cH,V_n)_{L_1[0,2]}=(n+1)d_n(\cH)_{L_1[0,2]},
\end{eqnarray*}
and the proof is completed.
\hfill $\Box$
\end{example}
It follows then from Lemma \ref{Albert} and Theorem \ref{TB} that the Lipschitz width of $\cH$ decays exponentially, while its Kolmogorov width decays like  $n^{-1}$.  While this is a good example, one may argue that this different behavior is due to the fact that $\cH$ is not convex. It is a well known fact that for every compact set $\cK$ we have 
$$
d_n(\cK)_X=d_n(\cK_c)_X, \quad \hbox{where} \quad \cK_c={\rm conv}(\cK\cup(-\cK))
$$
 is the minimal convex centrally symmetric set that contains $\cK$. Therefore, a more suitable example would be one when $\cK$ is a 
 convex, centrally symmetric set.  We discuss such case in  Example \ref{Example 2}.

\bigskip
\noindent
\begin{example}  \label{Example 2}  \rm
Consider the sequence $\sigma=(\sigma_j)_{j=1}^\infty$,  with $\sigma_j=(\log_2(j+1))^{-1/2}$, and the corresponding  linear   map 
 on sequences, $D_\sigma:\ell_1\to\ell_2$,  
$
\ell_1:=\{x=(x_1,x_2,\ldots):\,\sum_{j=1}^\infty|x_j|<\infty\},
$
defined as
$$
 D_\sigma (x)=y, \quad \hbox{where}\quad y_j=\sigma_j x_j, \quad j=1, 2, \ldots.
$$
Let us denote by $\cK_\sigma\subset \ell_2$ the image of the unit ball in   $\ell_1$ under this map, namely,
\be
\label{setK}
\cK_\sigma:=\{y\in\ell_2 :\ y_j =\sigma_j x_j,\, \mbox{ where}\ \sum_{j=1}^\infty |x_j|\leq 1\}=\{y\in \ell_2 \ :\ \sum_{j=1}^\infty |y_j|\sqrt{\log_2 ( j+1)} \leq 1\}.
\ee
The set $\cK_\sigma$ is a convex, centrally symmetric subset of $\ell_2$  for which
 $$
\left\{ \frac{\pm e_j}{\sqrt{\log_2(j+1)} } \right\}_{j=1}^\infty\subset \cK_\sigma.
$$
It follows from Proposition 3.1 in \cite{Kuhn} that
\be 
\label{entropyKuhn}
\e_n(\cK_\sigma)_{ \ell_2}\asymp n^{-1/2},\quad  n=1,2,\ldots,
\ee
which combined  Theorem \ref{TB} shows that 
 $d_n^{\gamma}(\cK_\sigma)_{ \ell_2}\leq Cn^{-1/2}$ with  $\gamma=2\rad(\cK_\sigma)=2$.
 On the other hand,    we show in the next lemma 
 that its Kolmogorov width $d_n(\cK_\sigma)_{ \ell_2}$ behaves as $d_n(\cK_\sigma)_{\ell_2}\asymp (\log_2 n)^{-1/2}$. 
 \begin{lemma}
 The Kolmogorov width of the compact set $\cK_\sigma$ defined in {\rm \eref{setK}} is
 $$
 d_n(\cK_\sigma)_{ \ell_2}\asymp (\log_2 n)^{-1/2}, \quad  n=2,3, \ldots.
 $$
 \end{lemma}
 \noindent
 {\bf Proof:} 
 Clearly, 
 $$
 d_n(\cK_\sigma)_{\ell_2}\leq
  \sup_{x\in \cK_\sigma}\dist(x,\span\{e_j\}_{j=1}^n)_{\ell_2} =\frac{1}{\sqrt{\log_2(n+2)}}.
$$
 To prove the  inequality from below, we  fix $\epsilon>0$ and denote by $X_n$ the  $n$ dimensional   subspace for which
 $$
 \sup_{x\in \cK_\sigma}\dist(x,X_n)_{ \ell_2}\leq (1+\epsilon)d_n(\cK_\sigma)_{ \ell_2}.
 $$
If  $P$ is the orthogonal projection onto $\ell_2^{2n}:=\span\{e_j\}_{j=1}^{2n}$ and 
$\widetilde X_n :=P(X_n)$, then
\be
\label{ff1}
d_n(P(\cK_\sigma))_{ \ell_2}\leq  \sup_{x\in P(\cK_\sigma)}\dist(x, \widetilde X_n)_{ \ell_2}\leq (1+\epsilon)d_n(\cK_\sigma)_{ \ell_2}.
\ee
Since  $\cP_n:=\frac{1}{\sqrt{\log_2(2n+1)}}\mathrm{conv}\{\pm e_j\}_{j=1}^{2n}\subset P(\cK_\sigma)$, we have
\be
\label{ff2}
d_n(\cP_n)_{ \ell_2}\leq d_n(P(\cK_\sigma))_{ \ell_2},
\ee
and from Stechkin's theorem  \cite[Ch. 13 Th.3.3] {LGM} we know that
\be
\label{ff3}
d_n(\cP_n)_{\ell_2}=  \frac{1}{\sqrt 2}\left(\log_2(2n+1)\right)^{-1/2}.
\ee
Combining \eref{ff1}, \eref{ff2} and \eref{ff3}  gives
$$
\frac{1}{\sqrt 2(1+\e)}\left(\log_2(2n+1)\right)^{-1/2}\leq  d_n(\cK_\sigma)_{\ell_2}.
$$
Since  $\epsilon>0$ is arbitrary, we obtain
$$
C(\log_2 n)^{-1/2}\leq  d_n(\cK_\sigma)_{\ell_2},
$$
which completes   the proof.
\hfill $\Box$
\end{example}

\section{Comparison between Lipschitz  and  stable manifold widths}
\label{S6}

 Let us recall the definition of   {\it manifold} width  $\delta_n(\cK)_X$  for the compact set $\cK\subset X$, see \cite{DHM,DKLT}, 
 \be
 \label{manwidth}
 \delta_n(\cK)_X:=\inf_{ a,M} \sup_{f\in\cK}\|f-M(a(f))\|_X,
 \ee
 where the infimum is taken over all  mappings $a:\cK\to \R^n$ and $M:\R^n\to X$ with  $a$ continuous on $\cK$ and $M$ continuous on $\R^n$.  A comparison between manifold widths and other types of nonlinear widths was given in \cite{DKLT}.    There is also another concept, called  {\it stable manifold }width $\delta_{n,\gamma}^*(\cK)_X$ of the compact set $\cK\subset X$, { see \cite{CDPW},}
   defined as
 \be
 \label{stablewidth}
 \delta_{n,\gamma}^*(\cK)_X:= \inf_{a,M,\|\cdot\|_{Y_n}} \sup_{f\in\cK}\|f-M(a(f))\|_X,
 \ee
 where now the infimum is taken  over all maps    $a:\cK\to (\R^n,\|\cdot\|_{Y_n})$, $M:(\R^n,\|\cdot\|_{Y_n})\to X$, and norms $\|\cdot\|_{Y_n}$ on $\R^n$, with  $a,M$ 
being  $\gamma$-Lipschitz.  We discuss in this section the relation between stable manifold widths and the Lipschitz widths. The next theorem shows that for any compact set $\cK\subset X$, the Lipschitz widths are smaller than the stable manifold widths.

\begin{theorem}
\label{TS}
For every compact set $\cK\subset X$, every $n\geq 1$, and every $\gamma>0$, we have
\be
\label{compStab}
d^{\gamma^2{\diam}(\cK)}_n(\cK)_X\leq \delta^*_{n,\gamma}(\cK)_X.
\ee
\end{theorem}
\noindent
{\bf Proof:}   We choose $\epsilon>0$, and let $a:\cK\to(\R^n,\|\cdot\|_{Y_n})$ and  $M:(\R^n,\|\cdot\|_{Y_n})\to X$ be two $\gamma$-Lipschitz 
mappings with respect to a norm $\|\cdot\|_{Y_n}$ in $\R^n$ such that
for every $f\in \cK$,
\be
\label{llk}
\|f-M\circ a(f)\|_X\leq \delta^*_{n,\gamma}(\cK)_X+\epsilon.
\ee
  For every $f_1,f_2\in \cK$ we have
$$
\|a(f_1)-a(f_2)\|_{Y_n}\leq \gamma\|f_1-f_2\|_X,
$$
which implies 
$$
\diam (\cA)\leq \gamma\diam(\cK), \quad \hbox{where}\quad \cA:=a(\cK).
$$	
We fix an element $f_0\in\cK$ and define the mapping $\Phi:(B_{Y_n},\|\cdot\|_{Y_n})\to X$ as
$$
\Phi(y):=M(a(f_0)+\gamma\diam(\cK)y), \quad a(f_0)\in \R^n.
$$
Note that 
$\Phi$ is a $\gamma^2\diam(\cK)$-Lipschitz mapping. For   each $f\in \cK$ we define
$$
 y(f):=
\frac{1}{\gamma\diam(\cK)}(a(f)-a(f_0))
\in B_{Y_n}.
$$
An easy calculation shows that  $\Phi( y(f))=M\circ a(f)$, so
 $$
 \|f-\Phi( y(f))\|_X=\|f-M\circ a(f)\|_X \leq \delta^*_{n,\gamma}(\cK)_X+\epsilon,
 $$
 where we have used \eref{llk}. 
 Therefore we obtain
 $$
 \sup_{f\in \cK}\inf_{y\in B_{Y_n}} \|f-\Phi(y)\|_X\leq \delta^*_{n,\gamma}(\cK)_X+\epsilon\quad \Rightarrow\quad 
  d_n^{\gamma^2\diam(\cK)}(\cK)_X\leq \delta^*_{n,\gamma}(\cK)_X+\epsilon.
 $$
 Since $\epsilon$ is arbitrary,  \eref{compStab} holds  and the proof is completed.
\hfill $\Box$

\begin{theorem}
For every Banach space $X$ and for every $n\in \N$ there exist   compact sets $\cK\subset X$ such that
for every $\gamma>0$, 
$$ 
\delta^*_{n,\gamma}(\cK)_X\geq \delta_n(\cK)_X\geq 1,\quad \hbox{while}\quad  \lim_{\gamma\to\infty} d^\gamma_n(\cK)_X=0.
$$
\end{theorem}
\noindent
{\bf Proof:} Let us  fix $n\in \N$ and consider the  compact set 
$$
\cK:=S_{n+1}\subset X_{n+1}\subset X,
$$
where $S_{n+1}$ is
  the boundary of the unit sphere of an $(n+1)$-dimensional subspace $(X_{n+1},\|\cdot\|_X)$ of $X$.  
By the Borsuk  theorem, see \cite{Bo,LGM},  we have that for any continuous map $a:\cK \to \R^n$, there exists 
$f_0\in \cK$ such that $a(f_0)=a(-f_0)$, and thus  for any map $M :\R^n \rightarrow X$ we have 
$
M(a(f_0))=M(a(-f_0)).
$
Then,  since $\|f_0\|_X=1$ and $f_0,-f_0\in \cK$, we have the inequality
\begin{eqnarray}
2&=&\|f_0-(-f_0)\|_{X}=\|f_0-M(a(f_0))+(M(a(-f_0))-(-f_0)\|_{X}\\ \nonumber
&\leq& \|f_0-M(a(f_0))\|_X+\|M(a(-f_0))-(-f_0)\|_{X}\leq 
 2\sup_{f\in \cK}\|f-M(a(f))\|_X
\end{eqnarray}
 for all  mappings $a:\cK\to \R^n$ and $M:\R^n\to X$ with  $a$ continuous on $\cK$ and $M$ continuous on $\R^n$.
So $\delta_n(\cK)_X\geq 1$, and therefore 
$\delta^*_{n,\gamma}(\cK)_X\geq \delta_n(\cK)_X\geq 1$ for any $\gamma>0$. 
 On the other hand, since $\cK$ is compact  (and thus totally bounded), we have that $\lim_{\gamma\to\infty} d^\gamma_n(\cK)_X=0$ because of Lemma \ref{Lem1}.
\hfill $\Box$

\section{Relation to neural networks}
\label{S7}

In this section, we discuss deep neural network approximation (DNNA) by feed-forward ${\rm ReLU}$ neural networks (NN) of constant width $W\geq 2$ and depth $n$, whose parameters have absolute values bounded by $1$. 
We will show  that  the approximation tools provided by these NNs   are in fact
Lipschitz mappings 
$$
\Phi:(B_{\ell_\infty^{\tilde n}},\|\cdot\|_{\ell_\infty^{\tilde n}})\to C(\Omega), \quad \Omega=[0,1]^d, \quad \tilde n=Cn, \quad C=C(W),
$$
with Lipschitz constant $\gamma_n=C'nW^n$. Therefore, a theoretical benchmark  for the  performance of the DNNA for a class  $\cK\subset X$
is given by the  Lipschitz width $d^{\gamma_n}_{\tilde n}(\cK)_X$.  This observation  motivates our investigation of Lipschitz widths whose Lipschitz constant depends on $n$.

\subsection{Deep neural networks as Lipschitz mappings}   
Let us  first recall that a DNNA  of a function $f\in C(\Omega)$,   $\Omega\subset \R^d$,  via feed-forward NN with activation function 
$\sigma:\R\to\R$, constant width $W$ and depth $n$ is in fact an approximation to $f$  by the family of functions  
$$
\Sigma_n:=\{\Phi(y): \,\,y\in \R^{\tilde n},\,\,\tilde n=\tilde n(W,n)=Cn\} \subset C(\Omega).
$$ 
For  each $y\in \R^{\tilde n}$, $\Phi(y) \in C(\Omega)$ is a continuous function  $\Phi(y):\Omega\to \R$ of the form
\be
\label{NN1}
\Phi(y):=A^{(n)}\circ\bar\sigma\circ A^{(n-1)}\circ\ldots\circ \bar\sigma\circ A^{(0)},
\ee
with  $A^{(0)}:\R^d\to\R^W$, $A^{(\ell)}:\R^{W}\to\R^{W}$, $\ell=1,\ldots,n-1$,  and $A^{(n)}:\R^W\to\R$ being affine mappings, and 
$\bar\sigma:\R^W\to\R^W$ given by 
$$
\bar\sigma(x_{j+1},\ldots,x_{j+W})=(\sigma(x_{j+1}),\ldots,\sigma(x_{j+W})).
$$
The argument $y$ of $\Phi$ is a vector in $\R^{\tilde n}$  that consists  of the entries of the matrices and offset   vectors (biases) of the affine mappings $A^{(\ell)}$,
$\ell=0,\ldots,n$.  We  order these entries  in such a way that the entries of $A^{(\ell)}$
 appear before those of $A^{(\ell+1)}$ and the ordering  for each $A^{(\ell)}$ is done in  the same way.
 Before going further, 
  we need to specify a norm $\|\cdot\|_{Y_{\tilde n}}$ to be used for $\R^{\tilde n}$.  
  We take this norm to be the $\ell_\infty^{\tilde n}:=\ell_\infty(\R^{\tilde n})$ norm, that is,   $\|y\|_{\ell_\infty^{\tilde n}}:=\max_{1\le i\le {\tilde n}}|y_i|$.  
  This choice is not optimal for obtaining the best constants in our estimates but it will simplify the exposition that follows. Also, when considering
  vector functions $g=(g_1,\dots,g_W)$  from $C(\Omega)$, we use the notation
   $$
  \|g\|:=\max_{1\le i\le W} \|g_i\|_{C(\Omega)}.
  $$

 It was proven in \cite{DHP} that 
  if $B$ is any finite ball in $\ell_\infty(\R^{\tilde n})$ and $\sigma(t)={\rm ReLU}(t)=\max\{t,0\} =t_+$, then 
$\Phi:B\to C(\Omega)$ is a $\gamma$-Lipschitz mapping with $\gamma$  depending only on $B,W,n$, and $d$. 
 In fact,  $\Phi$ is a $\gamma$-Lipschitz map on any bounded set. Here, we will investigate in detail the Lipschitz constant $\gamma$ in the case 
 when $B$ is the unit ball $(B_{\ell_\infty^{\tilde n}}, \|\cdot\|_{\ell_\infty^{\tilde n}})$. More precisely,  the following theorem holds.
  \begin{theorem}
  \label{NN}
  The mapping $\Phi:(B_{\ell_\infty^{\tilde n}}, \|\cdot\|_{\ell_\infty^{\tilde n}})\to C(\Omega)$, with $\Omega=[0,1]^d\subset \R^d$, defined in 
  {\rm \eref{NN1}} with $\sigma={\rm ReLU}$ is a $C'nW^n$-Lipschitz mapping, that is
  $$
  \|\Phi(y)-\Phi(y')\|_{C(\Omega)}\leq C'nW^n\|y-y'\|_{\ell_\infty^{\tilde n}}, \quad  y,y'\in B_{\ell_\infty^{\tilde n}},
  $$ 
  where $C'=C'(d)$ is a constant depending on $d$.
  \end{theorem}
  \noindent
 {\bf Proof:}  Let  $y$ and $y'$ be the entries of the affine mappings $A^{(j)}(\cdot):=A_j(\cdot)+b^{(j)}$,  $j=0,\ldots,n$,  and 
 $A'^{(j)}(\cdot):=A'_j(\cdot)+b'^{(j)}$,  $j=0,\ldots,n$, respectively, ordered in a predetermined way.
 We fix $x\in \Omega$ and  denote  by 
 $$
  \eta^{(0)}(x) := {\rm \overline{ReLU}} (A_0x+b^{(0)}),\quad  \eta'^{(0)}(x) := {\rm \overline{ReLU}} (A'_0x+b'^{(0)}),
$$
$$  
\eta^{(j)}:= {\rm \overline{ReLU}}(A_j\eta^{(j-1)}+b^{(j)}),\quad  \eta'^{(j)}:= {\rm \overline{ReLU}}(A'_j\eta'^{(j-1)}+b'^{(j)}), \quad j=1,\ldots,n-1,
$$
where $A_j,A_j',b^{(j)},b'^{(j)}$, $j=1,\ldots,n-1$,  are the respective $W\times W$ matrices and bias vectors,  associated to $y$ and $y'$,
and
$$
\eta^{(n)} := A_n\eta^{(n-1)}+b^{(n)},\quad  \eta'^{(n)}:=A'_n\eta'^{(n-1)}+b'^{(n)}.
$$
Note  that  since $\|y\|_{\ell_\infty^{\tilde n}}\leq 1$,
$$
\|\eta'^{(0)}\|\leq (d+1)\|y\|_{\ell_\infty^{\tilde n}}\leq d+1, \quad 
\|\eta'^{(j)}\|\leq
  (W\|\eta'^{(j-1)}\|+1)\|y\|_{\ell_\infty^{\tilde n}}\leq W\|\eta'^{(j-1)}\|+1, \quad j=1,\ldots,n.
  $$
One can show by induction that for $j=1,\ldots,n$, 
\be
\label{bound}
\|\eta'^{(j)}\|\leq W^{j}d+\sum_{k=0}^{j}W^k\leq (d+2)W^{j}.
\ee
Note that the above inequality also holds for $j=0$.
 Next, since  ReLU is a Lip 1 function, we have
 $$
 \|\eta^{(0)}(x)-\eta^{'(0)}(x)\|
 \leq \|(A_0-A_0')x\|+\|b^{(0)}-b'^{(0)}\|
 \leq (d+1)\|y-y'\|_{\ell_\infty^{\tilde n}}=:C_0\|y-y'\|_{\ell_\infty^{\tilde n}},
 $$
 and therefore $ \|\eta^{(0)}-\eta^{'(0)}\|\leq C_0\|y-y'\|_{\ell_\infty^{\tilde n}}$.
   Suppose we have proved that 
$$
 \|\eta^{(j-1)}-\eta'^{(j-1)}\|\leq C_{j-1}\|y-y'\|_{\ell_\infty^{\tilde n}}.
$$  
It follows that 
  \begin{eqnarray} 
  \label{rec1}
  \nonumber
  \|\eta^{(j)} (x)-\eta'^{(j)} (x)\|&\le& \|A_j\eta^{(j-1)} (x)+b^{(j)} -A_j'\eta'^{(j-1)} (x)-b'^{(j)}\|
  \\ \nonumber
  &\le&  \|A_j(\eta^{(j-1)} (x)-\eta'^{(j-1)} (x))\|+ \|(A_j-A'_j)\eta'^{(j-1)} (x)\|+\|b^{(j)}-b'^{(j)}\|
  \\ \nonumber
    &\le& 
W\|y\|_{\ell_\infty^{\tilde n}}\|\eta^{(j-1)}-\eta'^{(j-1)}\|+ W\|y-y'\|_{\ell_\infty^{\tilde n}}\|\eta'^{(j-1)}\|+\|y-y'\|_{\ell_\infty^{\tilde n}}\\
\nonumber
  &\le&
(WC_{j-1}+(d+2)W^{j}+1)\|y-y'\|_{\ell_\infty^{\tilde n}}
  \\ \nonumber
  &=:&
  C_{j}\|y-y'\|_{\ell_\infty^{\tilde n}},
  \end{eqnarray} 
  where  we  have used the induction hypothesis, the fact that $\|y\|_{\ell_\infty^{\tilde n}}\leq 1$, and the bound  \eref{bound} for $\|\eta'^{(j)}\|$. 
   Thus, we have obtained  that $\|\eta^{(j)}-\eta'^{(j)}\|\leq C_j\|y-y'\|_{\ell_\infty^{\tilde n}}$, and therefore the following   recursive relation
$$
C_{j}=WC_{j-1}+(d+2)W^{j}+1, \quad j=1,\ldots,n,
$$
 between the constants $C_j$, $j=1,\ldots,n$, where $C_0=d+1$.
 We then obtain
  $$
  C_{n}<(n+1)(d+2)W^{n}+\sum_{k=0}^{n-1}W^{k}<C'nW^{n}, \quad \hbox{with}\,\,C'=C'(d).
  $$
  Finally, we write 
 $$
\|\Phi(y)-\Phi(y')\|_{C(\Omega)}=\|\eta^{(n)}-\eta'^{(n)}\|\leq C_{n}\|y-y'\|_{\ell_\infty^{\tilde n}}<C'nW^{n}\|y-y'\|_{\ell_\infty^{\tilde n}},
  $$
  and the proof is completed.
\hfill $\Box$

We next discuss a Carl's type inequality  that is  similar to Lemma \ref{L3}, but is  for  the case when the Lipschitz constant $\gamma$ depends on $n$.

\begin{remark}
\label{RC}
If one follows the proof of Lemma  {\rm \ref{L3}} with the condition that $\gamma$ is not a constant, but $\gamma=\gamma_n=C'n^{\delta}\lambda^n$, where $C'>0$, $\delta\in \R$, and $\lambda>1$, one can show that
$$
d^{\gamma_n}_n(\cK)_X<  c_0 \frac{[\log_2n]^{\beta}}{n^{2\alpha}}, \quad \beta\in \R,\,\,\,\alpha>0\quad \Rightarrow\quad 
 \e_{m}(\cK)_X  < C\frac{[\log_2m]^{\beta}}{m^{\alpha}},\quad \hbox{where}\,\,m=cn^2,
$$
and $C,c$ are fixed constants, depending only on  $c_0$,  $\beta$, $\alpha$, $\delta$, $\lambda$ and $C'$.
Indeed, the proof follows from the fact that  for $\e=c_0[\log_2n]^{\beta}n^{-2\alpha}$ we have
$$
N_{2\e}(\cK)\leq  \left(\frac{3\gamma}{\e}\right)^n=
(3C'c_0^{-1}\lambda^n [\log_2 n]^{-\beta}n^{\delta+2\alpha})^n <  
2^{cn^2}, 
$$
and therefore
$$
\e_{cn^2}(\cK)_X<   2c_0 [\log_2 n]^{\beta}n^{-2\alpha}.
$$
Setting $m=cn^2$   i.e. $n=\sqrt{m/c}$ gives

$$
\e_m(\cK)_X\leq 2 c_0 [\log_2 \sqrt{m/c}]^\beta (m/c)^{-\alpha}= 2 c_0c^\alpha2^{-\beta}\frac{ [\log_2 m -\log_2 c]^\beta }{m^\alpha}<C^{\prime \prime}\frac{ [\log_2 m]^\beta }{m^\alpha},
$$
 which is what we wanted to show.
\end{remark}

\subsection{Lower bound for $d_n^{\gamma_n}(\cK)_X$}
Now that we know that DNNA is an approximation to a function $f$ by a  particular $\gamma_n$-Lipschitz mapping
with $\gamma_n=C'nW^n$, we can ask the question what are the limits of such  approximation. This question is answered by providing a lower bound for the Lipschitz width $d_n^{\gamma_n}(\cK)_X$ via the  next theorem  which is a modification of Theorem \ref{widthsfrombelow}. 

\begin{theorem}
\label{widthsfrombelownew}  
For any  compact set $\cK\subset X$ we consider the Lipschitz width $d_n^{\gamma_n}(\cK)_X$ with 
$\gamma_n=C'n^{\delta}\lambda^n$, $\delta\in\R$,  $\lambda >1$ and $C'>0$ being  fixed constants. Then the following holds:
\begin{enumerate}
\item   if for some constants $c_1>0,\alpha>0$ and $\beta\in \R$ we have 
$$
\e_n(\cK)_X> c_1 \frac{(\log_2 n)^\beta}{n^{\alpha}}, \quad n=1,2,\dots,
$$
 then there exists a constant $C>0$ such that
\be 
\label{widths(i)new}
d_n^{\gamma_n}(\cK)_X\geq C\frac{(\log_2 n)^{\beta}}{n^{2\alpha}}, \quad n=1,2,\dots.
\ee

\item  if for some constants $c_1>0,\alpha>0$  we have  
$$
\e_n(\cK) _X> c_1 (\log_2 n)^{-\alpha}, \quad n=1,2,\dots,
$$ 
then there exists a constant $C>0$ such that 
\be 
\label{widths(ii)new}
d_n^{\gamma_n}(\cK)_X\geq C(\log_2 n)^{-\alpha}, \quad n=1,2,\dots .
\ee

\end{enumerate}
\end{theorem}
\noindent
{\bf Proof:} 
We prove  the theorem  by  contradiction. We first concentrate on the proof of (i).
 If \eref{widths(i)new} does not hold for some constant $C$,  then
there exists a strictly  increasing sequence of integers $(n_k)_{k=1}^\infty$, such that
$$ 
p_k:= \frac{d_{n_k}^{\gamma_{n_k}}(\cK)n_k^{2\alpha}}{(\log_2 n_k)^{\beta}}\to 0\quad\hbox{as}\quad k\to \infty.
$$
Thus we can write
\be 
\nonumber
d_{n_k}^{\gamma_{n_k}} (\cK)  = \frac{p_k\left[\log_2 n_k\right]^{\beta}}{n_k^{2\alpha}}<
\frac{2p_k\left[\log_2 n_k\right]^{\beta}}{n_k^{2\alpha}}=:\delta_k, \quad  k=1,2,\dots.
\ee
Now  we apply Proposition \ref{carl2}  with $\eta_n=c_1(\log_2n)^\beta n^{-\alpha}$
and obtain
$$
c_1 \left[\log_2 ({ n_k} \log_2(3{\gamma_{n_k}}\delta_{k}^{-1}))\right ]^\beta  {n_k}^{-\alpha} \left[\log_2(3{\gamma_{n_k}}\delta_{k}^{-1})\right]^{-\alpha } \leq 4 \frac{
p_k\left[\log_2 {n_k}\right]^{\beta}}{ {n_k}^{2\alpha}},
$$
which we rewrite as
\be
\label{f1new}
p_k^{-1}\left[\log_2  {n_k} + \log_2 \log_2(3{\gamma_{n_k}}\delta_{k}^{-1})\right]^\beta  
\left[\log_2(3{\gamma_{n_k}}\delta_{k}^{-1})\right]^{-\alpha } \leq C_1  \left[\log_2 {n_k}\right]^{\beta}n_k^{-\alpha},\quad C_1=4/c_1.
\ee
Observe that
\begin{eqnarray*}  
\log_2(3{\gamma_{n_k}}\delta_{k}^{-1})&=&
 \log_2 {(1.5{\gamma_{n_k}})} +\log_2 (p_k^{-1} )+2\alpha \log_2 {n_k} -\beta \log_2 (\log_2 {n_k})\\
&=&
 \log_2 {(1.5C'n_k^{\delta}\lambda^{n_k})} +\log_2 (p_k^{-1}) +2\alpha \log_2 {n_k} -\beta\log_2 (\log_2 {n_k})\
,
\end{eqnarray*}
and therefore  for $k$ big enough we obtain
\be
\label{gqnew}
 \log_2 (3\gamma_{n_k}\delta_{k}^{-1})\leq 2\left[\log_2 (p_k^{-1}) +A  n_k\right ].
\ee
The latter inequality and  \eref{f1new} give
$$
 p_k ^{-1}\left[\log_2  {n_k} + \log_2 (\log_2(3\gamma_{n_k}\delta_{k}^{-1}))\right]^\beta \left[\log_2 (p_k^{-1}) +A n_k\right ]^{-\alpha}
\leq 2^{\alpha}C_1\left[\log_2 {n_k}\right]^{\beta}n_k^{-\alpha},
$$
which is equivalent to
\be
\label{f2new}
p_k^{-1}\left[\log_2  {n_k} + \log_2 (\log_2(3\gamma_{n_k}\delta_{k}^{-1}))\right]^\beta
\leq 2^\alpha C_1\left[\frac{\log_2{  (p_k^{-1})}}{n_k}+A\right ]^{\alpha}\left[\log_2 {n_k}\right]^{\beta}.
\ee

{\bf Case 1:}  $\beta\geq 0$.

\noindent
Note that since $\delta_k\to 0$ and $\gamma_{n_k}\to \infty$,  for $k$ big enough we have $\log_2( \log_2(3\gamma_{n_k}\delta_{k}^{-1}))>0$. Since   $\beta\geq 0$, we have 
$$
\left[\log_2  {n_k}\right]^\beta\leq \left[\log_2  {n_k} + \log_2( \log_2(3\gamma_{n_k}\delta_{k}^{-1}))\right]^\beta, 
$$
and therefore it follows from \eref{f2new} that
$$
p_k^{-1}\leq  2^\alpha C_1\left[\frac{\log_2 (p_k^{-1})}{n_k }+A \right ]^{\alpha}<C[\log_2 (p_k^{-1})]^{\alpha},
$$
which contradicts the fact that $f_{k}$ tends to zero (and thus $p_k^{-1}\to \infty$).

{\bf Case 2:}   $\beta< 0$. 

\noindent
In this case we rewrite \eref{f2new} and use \eref{gqnew} to obtain 
\begin{eqnarray*}
p_k^{-1}
&\leq &2^\alpha C_1\left[\frac{\log_2 (p_k^{-1})}{n_k}+A\right ]^{\alpha}
\left[1+ \frac{\log_2 (\log_2(3\gamma_{n_k}\delta_{k}^{-1}))}{\log_2  {n_k} }\right]^{-\beta}\\
&\leq&
 2^\alpha  C_1\left[\frac{\log_2 (p_k^{-1})}{n_k}+A\right ]^{\alpha}
\left[1+ \frac{\log_2 (2An_k+ 2 \log_2(p_k^{-1}))}{\log_2  {n_k} }\right]^{-\beta}.
\end{eqnarray*}
Next, we consider the following 2 cases.

{\bf Case 2.1:}   If for infinitely many values of $k$ we have $p_k^{-1}\leq cn_k$, then the above  inequality becomes
$$
p_k^{-1}\leq C,
$$
which contradicts with the fact that $p_k^{-1}\to \infty$ as $k\to \infty$.

{\bf Case 2.2:}   If for infinitely many values of $k$ we have $p_k^{-1}\geq cn_k$, then the above  inequality becomes
$$
p_k^{-1}\leq C\left[\log_2 (p_k^{-1})\right]^{\alpha}\left[\log_2 (\log_2(p_k^{-1}))\right]^{-\beta},
$$
which also contradicts with the fact that $p_k^{-1}\to \infty$ as $k\to \infty$.

To prove (ii) we repeat the argument for (i), namely, we assume that (ii) does not hold. Therefore there exists a strictly  increasing sequence of integers 
$(n_k)_{k=1}^\infty$, such that
$$ 
e_k:= d_{n_k}^{\gamma_n}(\cK)[\log_2 n_k]^{\alpha}\to 0\quad\hbox{as}\quad k\to \infty.
$$
We write
\be 
\label{widths(iii)new}
d_{n_k}^{ \gamma_{n_k}} (\cK)= e_k[\log_2 n_k]^{-\alpha}<
2e_k[\log_2 n_k]^{-\alpha}=:\delta_k, \quad k=1,2,\dots,
\ee
and use  Proposition \ref{carl2}  with $\eta_n=c_1(\log_2n)^{-\alpha}$
to derive
$$
c_1 \left[\log_2 (n_k \log_2(3\gamma_{n_k}\delta_{k}^{-1}))\right ]^{-\alpha}
 \leq 4e_k[\log_2 {n_k}]^{-\alpha}.
$$
The latter inequality is equivalent to
\be
\label{f11}
e_k^{-1}\leq C_1\left[1+ \frac{\log_2 (\log_2(3\gamma_{n_k}\delta_{k}^{-1}))}{\log_2  {n_k} }\right]^{\alpha } \leq
C_1\left[1+ c\frac{\log_2 (n_k+\log_2(e_k^{-1}))}{\log_2  {n_k} }\right]^{\alpha },\quad  C_1=4/c_1,
\ee
where we have used inequality similar to \eref{gqnew}.

{\bf Case 1:}   If for infinitely many values of $k$ we have $e_k^{-1}\leq cn_k$, then the above  inequality becomes
$$
e_k^{-1}\leq C,
$$
which contradicts with the fact that $e_k^{-1}\to \infty$ as $k\to \infty$.

{\bf Case 2:}   If for infinitely many values of $k$ we have $e_k^{-1}\geq cn_k$, then the above  inequality becomes
$$
e_k^{-1}\leq C\left[\log_2 (e_k^{-1})\right]^{\alpha},
$$
which also contradicts with the fact that $e_k^{-1}\to \infty$ as $k\to \infty$.
\hfill $\Box$

\subsection{Summary}

In this section we summarize our results for the Lipschitz widths $d_n^{\gamma_n}(\cK)_X$ and give several examples.
The following corollary holds.

\begin{cor}
\label{cNN}
Let  $\cK\subset X$ be a  compact subset of a Banach space $X$,  $n\in \N$, and $d_n^{\gamma_n}(\cK)_X$ be the   Lipschitz width for $\cK$ with Lipschitz constant  
$\gamma_n=C'n^{\delta}\lambda^n$, where $\delta\in \R$, $C'>0$  and $\lambda> 2$. 
\begin{enumerate}
\item For  $\alpha>0$, $\beta\in \R$, we have 
\be
\label{c1}
\e_n(\cK)_X\asymp \frac{[\log_2n]^\beta}{n^\alpha} \quad \Rightarrow\quad 
d_n^{\gamma_n}(\cK)_X\asymp \frac{[\log_2n]^\beta}{n^{2\alpha}};
\ee
\item For  $\alpha>0$,  we have 
\be
\label{c2}
\e_n(\cK)_X\asymp \frac{1}{[\log_2 n]^\alpha} \quad \Rightarrow\quad 
d_n^{\gamma_n}(\cK)_X\asymp \frac{1}{[\log_2 n]^{\alpha}}.
\ee
\end{enumerate}
\end{cor}
\noindent
{\bf Proof:} We first prove (i). Let us assume that 
$$
\e_n(\cK)_X\leq  C\frac{[\log_2n]^\beta}{n^\alpha},  
$$
holds. After using  \eref{power} from Theorem \ref{TB}, we obtain
\be
\label{aa}
d_n^{2^n\rad(\cK)}(\cK)_X\leq 2^{\beta}C\frac{[\log_2n]^\beta}{n^{2\alpha}}.
\ee
We now fix $n_0$ such that $C'n^{\delta}\lambda^n\geq 2^n\rad(\cK)$ for all $n\geq n_0$ 
(recall that $\lambda>2$). We apply
Lemma \ref{LL1}, (iv) to derive
\be
\label{aa1}
d_n^{\gamma_n}(\cK)_X\leq d_n^{2^n\rad(\cK)}(\cK)_X, \quad n\geq n_0.
\ee
Finally, it follows from \eref{aa} and \eref{aa1} that 
$$
d_n^{\gamma_n}(\cK)_X\leq C'\frac{[\log_2n]^\beta}{n^{2\alpha}}, \quad \hbox{for all}\,\, n,
$$
 provided the constant  $C'$ is chosen appropriately.
The other direction in \eref{c1} is the statement of Theorem \ref{widthsfrombelownew}, (i).
The proof of (ii) is similar and we omit it.
\hfill $\Box$

 Corollary \ref{cNN} provides a tool for giving lower bounds on how well a compact set (model class) $\cK$ can be approximated by a DNN. So far, one way to give such lower bounds is via VC dimension, see \cite{DHP}, \S 5.9, and the references therein,  which  
is  restricted to the case when approximation error is measured in the norm
  $\|\cdot\|_{C(\Omega)}$. Note that  Corollary \ref{cNN} can be applied in the case of $L_p$ approximation 
  when $p\neq \infty$. For example, if  $B^s_ q(L_\tau(\Omega))$,  $\Omega=[0,1]^d$, is any Besov space that lies above the Sobolev embedding line for $L_p(\Omega)$, then it is proven in \cite{DS} that
  $$
\e_n(U(B^s_q(L_\tau(\Omega))))_{L_p(\Omega)}\asymp n^{-s/d},   
  $$
  where $U(B^s_q(L_\tau(\Omega)))$ is the unit ball of $B^s_q(L_\tau(\Omega))$. Then, according to Theorem \ref{NN} and Corollary \ref{cNN}, we have 
  $$
 \dist(U(B^s_q(L_\tau(\Omega))), \Sigma_n)_{L_p(\Omega)}\geq d^{\gamma_n}_{\tilde n}(U(B^s_q(L_\tau(\Omega))))_{L_p(\Omega)}\geq Cn^{-2s/d}.
  $$
 In particular, we recover the estimate, see (5.18) in \cite{DHP}
 $$
\dist(U(B^s_q(L_\infty(\Omega))), \Sigma_n)_{C(\Omega)}\geq Cn^{-2s/d}.
 $$
  Note that, in contrast to stable manifold widths,  the Lipschitz  widths do not shed a light on  the numerical aspect of 
 this approximation, that is, they do  not give  even a theoretical algorithm of how to design the approximant.

\vskip .1in
\noindent
{\bf Affiliations:}

\noindent
Guergana Petrova, Department of Mathematics, Texas A$\&$M University, College Station, TX 77843,  gpetrova$@$math.tamu.edu
\vskip .1in
\noindent
Przemys{\l}aw Wojtaszczyk, Institut of Mathematics Polish Academy of Sciences, ul. 
{\'S}niadeckich 8,  00-656 Warszawa, Poland, wojtaszczyk$@$impan.pl

\end{document}